\input amstex.tex
\documentstyle{amsppt}
\UseAMSsymbols\TagsAsMath\widestnumber\key{ASSSS}
\magnification=\magstephalf\pagewidth{6.2in}\vsize8.0in
\parindent=6mm\parskip=3pt\baselineskip=16pt\tolerance=10000\hbadness=500
\NoRunningHeads\loadbold\NoBlackBoxes\nologo
\def\today{\ifcase\month\or January\or February\or March\or April\or May\or June\or
     July\or August\or September\or October\or November\or December\fi \space\number\day,
     \number\year}

\def\ns#1#2{\mskip -#1.#2mu}\def\operator#1{\ns07\operatorname{#1}\ns12}
\def\pa{\partial}\def\na{\nabla\!}
\def\dt{{\ns18{\Cal D}_{\ns23 t\!} }}\def\pdt{{\ns18\partial\ns18{\Cal D}_{\ns23 t\!} }}

\def\div{\operator{div}}
\def\dist{\operatorname{dist} }

\topmatter
\title
Well-posedness for the linearized motion of a compressible liquid
with free surface boundary.
\endtitle
\author  Hans Lindblad
\endauthor
\thanks  The author was supported in part by the National
Science Foundation.  \endthanks
\address
University of California at San Diego
\endaddress
\date \today\enddate
\email
lindblad\@math.ucsd.edu
\endemail
\endtopmatter

\document

\head{1. Introduction}\endhead
We consider Euler's equations
$$
\rho \big(\partial_t +V^k\partial_{k} \big) v_j=-\partial_{j\,} p,\quad j=1,...,n
\quad \text{in}\quad {\Cal D},\qquad \text{where}\quad
\partial_i =\partial/\partial x^i,\tag 1.1
$$
describing the motion of a perfect compressible fluid body in vacuum:
$$
 (\partial_t+V^k\partial_{k}) \rho+ \rho \div V=0 ,\qquad \div V=\partial_{k} V^k
\qquad \text{in}\quad {\Cal D},\tag 1.2
$$
where  $V^{k \!}=\delta^{ki} v_i=v_k$ and we use the summation
convention over repeated upper and lower indices. Here the
velocity $V\!=(V^1\!,...,V^n)$, the density $\rho$ and the domain
${\Cal D}\!=\cup_{0\leq t\leq T} \{t\}\!\times {\Cal D}_t$, ${\Cal
D}_t\!\subset\Bbb R^n$ are to be determined. The pressure
$p=p\,(\rho)$ is assumed to be a given strictly increasing smooth
function of the density. The boundary $\pa{\Cal D}_t$ moves with
the velocity of the fluid particles at the boundary. The fluid
body moves in vacuum so the pressure vanishes in the exterior and
hence on the boundary. We therefore also require the boundary
conditions on 
$\pa {\Cal D}=\cup_{0\leq t\leq T} \{t\}\!\times \pa {\Cal D}_t$: 
$$
\align   (\partial_t+V^k\partial_{k})|_{\partial {\Cal D}}&\in T(\partial {\Cal
D}),\tag 1.3\\
p=0,\quad & \text{on}\quad \partial {\Cal D}.\tag 1.4
\endalign
$$

Constant pressure on the boundary leads to energy conservation and it is
 needed for the linearized equations
to be well posed. Since the pressure is assumed to be a strictly
increasing function of the density we can alternatively think of
the density as a function of the pressure and for physical reasons
this function has to be non negative. Therefore
 the density has to be a non negative constant
$\overline{\rho}_0$ on the boundary and we will in fact assume
that $\overline{\rho}_0\!>\!0$, which is the case of liquid. 
We hence assume that 
$$
p(\overline{\rho}_0)=0\qquad\text{and}\qquad   p^\prime(\rho)>0,
\qquad\text{for}\quad \rho\geq \overline{\rho}_0,
\qquad\text{where}\quad \overline{\rho}_0>0\tag 1.5 
$$
From a physical point of view one can alternatively think of
the pressure as a small positive constant on the boundary. 
By thinking of the density as function of the pressure 
the incompressible case can be thought of 
as the special case of constant density function.

The motion of the surface of the ocean is described by the above
model. Free boundary problems for compressible fluids are also
of fundamental importance in astrophysics since they describe
stars. The model also describes the case of one fluid surrounded
by and moving inside another fluid.  For large massive bodies like
stars gravity helps holding it together and for
smaller bodies like water drops surface tension helps
holding it together. Here we neglect the influence of gravity
which will just contribute with a lower order term and we neglect
surface tension which has a regularizing effect.

Given a bounded domain
${\Cal D}_0\subset \Bbb R^n$, that is homeomorphic to the unit ball,
and initial data $V_0$ and $\rho_0$,
we want to find a set
${\Cal D}\subset [0,T]\times \Bbb R^n$,
a vector field $V$ and a function $\rho$,
solving (1.1)-(1.4) and satisfying the initial conditions
$$
\align \{x;\, (0,x)\in {\Cal D}\}&={\Cal  D}_0,\tag 1.6\\
V=V_0,\quad\rho=\rho_0\quad &\text{on}\quad \{0\}\times {\Cal D}_0.\tag 1.7
\endalign
$$
In order for the initial-boundary value problem (1.1)-(1.7) to be
solvable initial data (1.7) has to satisfy certain compatibility
conditions at the boundary. By (1.2), (1.4) also implies that
$\div V\big|_{\pa{\Cal D}}=0$. We must therefore have
$\rho_0\big|_{\pa{\Cal D}_0}=\overline{\rho}_0$ and $\div
V_0\big|_{\pa{\Cal D}_0}=0$.
Furthermore, taking the divergence of (1.1) gives an equation for
$(\pa_t+V^k\pa_k )\div V$ in terms of only space derivatives of
$V$ and $\rho$, which leads to further compatibility conditions.
In general we say that initial data satisfy the compatibility
condition of order $m$ if there is a formal power series solution
in $t$, of (1.1)-(1.7) $(\tilde{\rho},\tilde{V})$, satisfying
$$
(\pa_t +\tilde{V}^k\pa_k)^j (\tilde{\rho}-\overline{\rho}_0)\big|_{\{0\}\times
\pa{\Cal D}_0}=0,\qquad\quad j=0,..,m-1\tag 1.8
$$

Let ${\Cal N}$ be the exterior unit normal to the
free surface $\pa{\Cal D}_t$.
Christodoulou\cite{C2} conjectured the initial value problem (1.1)-(1.8),
is well posed in  Sobolev spaces under the assumption
$$
\na_{\Cal N}\, p\leq -c_0 <0,\quad\text{on}\quad \partial {\Cal D},
\qquad \text{where}\quad \na_{\Cal N}={\Cal N}^i\partial_{x^i}. \tag 1.9
$$
Condition (1.9) is a natural physical condition.
It says that the pressure and hence the density is larger
in the interior than at the boundary.
Since we have assumed that the pressure vanishes or is close to zero
at the boundary this is therefore related to the fact that the pressure of a
fluid has to be positive.

In general it is possible to prove local existence for analytic data
for the free interface between two fluids. 
However, this type of problem might be subject to instability in Sobolev norms, 
in particular Rayleigh-Taylor instability, which occurs when a heavier fluid is
on top of a lighter fluid. 
Condition (1.9) prevents Rayleigh-Taylor instability from occurring. 
Indeed,  if this condition is violated  Rayleigh-Taylor instability occurs
in a linearized analysis.  

In the irrotational incompressible case the physical condition
(1.9) always hold, see \cite{W1,2,CL}, and \cite{W1,2} proved
local existence in Sobolev spaces in that case. \cite{W1,2}
studied the classical water wave problem describing the motion of
the surface of the ocean and showed that the water wave is not
unstable when it turns over. Ebin\cite{E1} showed that
the general incompressible problem is ill posed in Sobolev spaces when the
pressure is negative in the interior and the physical condition is
not satisfied. Ebin\cite{E2} also announced a local
existence result for the incompressible problem with
surface tension on the boundary which has a regularizing effect
so (1.9) is not needed then. 

In \cite{CL}, together with Christodoulou, we proved {\it
a priori} bounds in Sobolev spaces in the general incompressible
case of non vanishing curl, assuming the physical condition
(1.9) for the pressure. We also showed that the Sobolev
norms remain bounded as long as the physical condition hold and
the second fundamental form of the 
free surface and the first order derivatives of the velocity are bounded.
Usually, existence follows from similar bounds
for some iteration scheme, but the bounds in \cite{CL}
used all the symmetries of the equation 
and so only hold for modifications that preserve all the symmetries. 
In \cite{L1} we showed existence for
the  linearized equations and in \cite{L2} we proved
local existence for the nonlinear incompressible problem with non
vanishing curl, assuming that (1.9) holds initially.

For the corresponding compressible free boundary problem with
non-vanishing density on the boundary, there are however in
general no previous existence or well-posedness results.
Relativistic versions of these problems have been studied in
\cite{C1,DN,F,FN,R} but solved only in special cases. The methods
used for the irrotational incompressible case use that the
components of the velocity are harmonic to reduce the equations to
equations on the boundary and this does not work in the
compressible case since the divergence is non vanishing and the
pressure satisfies a wave equation in the interior. To be able to
deal with the compressible case one therefore needs to use
interior estimates as in \cite{CL,L1}. 
Let us also point out that in nature one expects fluids to be compressible,
e.g.${}_{\!\!}$  water satisfies (1.5), see \cite{CF}. 
For the general relativistic equations there is no special case 
corresponding to the incompressible case.

Here we show existence for the linearized equations and estimates for these
in Sobolev spaces in the general compressible case (1.1)-(1.7), assuming that
(1.8) and (1.9) hold. This can be considered as a linearized stability result,
showing that small perturbations of initial conditions in Sobolev spaces leads to
small perturbations for finite times. 
Furthermore, in a forthcoming paper \cite{L3} we use existence 
and estimates for the inverse of the linearized operator to
prove existence for the nonlinear problem using the Nash-Moser technique
also in the compressible case.
The existence proof here uses the orthogonal decomposition of a vector
field into a divergence free part and a gradient of a function that vanishes
on the boundary.
For the divergence free part we get an equation of the type studied
in \cite{L1} and for the divergence we get a wave equation
on a bounded domain with Dirichlet boundary conditions.
The interaction terms between these equations are lower order
so if we put up an iteration, the equations decouple for the
the new iterate and the previous iterates only enter in the lower order terms.

Existence of solutions for the wave equation on a bounded domain is well known.
However dealing with the divergence free part of the equation requires
the techniques developed in \cite{L1}. Here we use a generalization of the
existence theorem in \cite{L1} to the case when the divergence of the solution
we linearize around is non vanishing.
In \cite{L1} we showed that the linearized incompressible Euler's  equations
becomes an evolution equation for what we called the normal operator.
The normal operator is unbounded and not elliptic
in the case of non vanishing curl. It is however positive assuming the physical
condition (1.9) and this leads to existence. Up to lower order terms,
the projection of the linearized compressible Euler's equations onto
divergence free vector fields becomes the linearized incompressible Euler's
equations.

As pointed out above, the positivity of the pressure (1.9) leads to the 
positivity
of the normal operator, introduced in \cite{L1}. It appears that this condition
is needed for the well-posedness also in the compressible case since the
divergence free part essentially decouples from the divergence.
In fact, the compressible case was the
main motivation for formulating (1.9) since in
that case it is clear that the pressure has to be positive
and in nature one expects fluids to be slightly compressible.

In order to formulate the linearized equations one has to parametrize the boundary.
Let us therefore express Euler's equations in the
 Lagrangian coordinates given by following the flow lines of the velocity vector 
field of the fluid particles. In these coordinates the boundary becomes fixed.
Given a domain ${\Cal D}_0$ in $\bold{R}^n$,
that is diffeomorphic to the unit ball $\Omega$, 
we can by a theorem in \cite{DM} find a
diffeomorphism $f_0:\Omega\to{\Cal D}_0$ with prescribed volume form 
$\det{(\pa f_0/\pa y)}$ up to a constant factor. 
Let ${\Cal D}$ and $v\!\in \! C({\Cal D})$ satisfy (1.3).
The Lagrangian coordinates $y$ are given by solving for the Eulerian coordinates
$x\!=\!x(t,y)\!=\!f_t(y)$
$$
\frac{d x(t,y)}{dt}=V(t,x(t,y)), \qquad x(0,y)=f_0(y),\qquad y \in \Omega.\tag 1.10
$$
Then $f_t:\Omega\to {\Cal D}_t$ is a diffeomorphism
and the boundary becomes fixed in the new $y$ coordinates.
Let
$$
 D_t=\frac{\partial }{\partial t}\Big|_{x=const}\!\!
+\,  V^k\frac{\partial}{\partial x^k}
=\frac{\partial }{\partial t}\Big|_{y=const}\!\!
\qquad\text{and}\qquad \partial_i=\frac{\partial}{\partial {x^i}}
=\frac{\partial y^a}{\partial x^i}
\frac{\partial}{\partial y^a}, \tag 1.11
$$
be the material derivative and partial differential operators expressed in the
Lagrangian coordinates.

In these coordinates  Euler's equation (1.1), the continuity equation (1.2) and
the boundary condition (1.4)  become:
$$
 D_t^2 x^i+\pa_{i\,}  h=0, \qquad
 D_t\, \rho +\rho\, \div V=0, \qquad \rho\big|_{\pa\Omega}=\overline{\rho}_0,\tag 1.12
$$
where the enthalpy
 $h=h(\rho)=\int_{\overline{\rho}_0}^\rho p^\prime(\rho)\rho^{-1}\, d\rho$
is a strictly increasing
function of $\rho$ ,  and $x$, $V=D_t\, x$
and $\rho$ are functions of $(t,y)\in[0,T]\times\Omega$. Furthermore,
$\rho$ can be determined from $x$:
$$
\rho=k\kappa^{-1},\qquad\text{where}\qquad \kappa=\det{(\pa x/\pa y)}
\qquad\text{and}\qquad k=\rho\kappa\big|_{t=0},\tag 1.13
$$
since $D_t\,\kappa =\kappa \div V$.
In (1.10) there is a choice of mapping $f_0$ and domain $\Omega$.
By \cite{DM} one can find a diffeomorphism with prescribed
volume form up to a constant between any  two diffeomorphic sets.
We therefore choose $\Omega$ to be the unit ball and
$\det{(\pa f_0/\pa y)}$ so $\det{(\pa f_0/\pa y)}\rho_0\!=\!k$ is any given fixed function
$k(y)$ that we take to be constant. Making this choice, initial data for $\rho$
is  part of the initial data for $x$.

The free boundary problem for Euler's equations (1.1)-(1.7), hence become an equation 
for $x(t,y)$: 
$$
 D_t^2 x^i+\pa_{i\,}  h=0, \qquad
 \rho=k\kappa^{-1}\!\!\!\!, \qquad \text{in}\quad [0,T]\times\Omega,\qquad
\kappa\big|_{\pa\Omega}=1,\quad\text{where}\quad 
\pa_i=\frac{\pa y^a}{\pa x^i}\frac{\pa}{\pa y^a},\tag 1.14
$$
$h\!=\!h({}_{\!}\rho{}_{\!})$ is a strictly increasing function of $\rho\!$
and $\rho\!=\!{}_{\!}\rho(_{\!}\kappa)$ is a function of 
$\kappa\!=\!\det_{\!}{(\pa x_{\!}/\pa y)}.\!$
Initial data are
$$
x\big|_{t=0}=f_0,\qquad D_t x\big|_{t=0}=V_0\tag 1.15
$$
In order for (1.14) to be solvable, initial data has to satisfy the constraints;
$\det{(\pa f_0/\pa y)}\big|_{\pa\Omega}=1$,
$\div V_0\big|_{\pa\Omega}=0$ and taking the divergence of (1.14)
gives an equation for $D_t \div V$ in terms of space
derivative of $x$ and $V=D_{t\,} x$ which leads to further conditions.
Since (1.14) gives $D_t^2 x$ in terms of space derivatives of $x$
we can obtain a formal power series solution in time $t$, $\tilde{x}$, to the first two
equations in (1.14) satisfying the initial conditions (1.15). The compatibility
condition of order $m$ is the requirement
that the formal power series solution up to terms of order $m$
satisfy the boundary condition in (1.14):
$$
D_t^j \big(\det{(\pa \tilde{x}/\pa y)}-1\big)\big|_{0\times\pa\Omega}=0,
\qquad j=0,...,m-1\tag 1.16
$$

Let us now derive the linearized equations.
(1.14) can be thought of as an equation $\Phi(x)=0$,
where $\Phi$ is a functional of $x(t,y)$ given by 
$\Phi(x)_i=D_t^2 x^i+\pa_i h$, for $1\leq i\leq n$,
where $h$ is a given function of 
$\kappa=\det{(\pa x/\pa y)}$ and $\pa_i$ are the differential operators in (1.14)
with coefficients depending on derivatives of $x$ as well, and 
$\Phi(x)_{n+1}=(\kappa-1)\big|_{\pa\Omega}$.
We assume that $x(t,y)$ is a given smooth solution of (1.14), i.e. $\Phi(x)=0$. 
Let $\overline{x}(t,y,r)$ be a smooth function
also of a parameter $r$, such that  $\overline{x}\big|_{r=0}=x$ and set
$\delta x=\pa \overline{x}/\pa r\big|_{r=0}$.
Then the linearized equations are the requirement on $\delta x$,
that $\overline{x}$
satisfies the equations (1.14) up to terms bounded by $r^2$ as $r\to 0$,
i.e. $\Phi^\prime(x)(\delta x)=\pa \Phi(\overline{x})/\pa r\big|_{r=0}=0$. 
If we replace $x$ in (1.14) by $\overline{x}$ and
apply $\delta =\pa /\pa r\big|_{r=0}$  we hence obtain the linearized
equations:
$$
 D_t^2\delta x^i +(\pa_i\pa_k h)\delta x^k
- \pa_i\big((\pa_k h)\delta x^k-\delta h\big)=0,
\qquad \delta h=-h^\prime(\rho)\rho\, \div\delta x,
\qquad \div\delta x\big|_{\pa\Omega}=0.\tag 1.17
$$
Here we used that $[\delta,\pa_i]=-(\pa_i\delta x^k)\pa_k$ and
$\delta \rho=-\rho \kappa^{-1}\delta\kappa=-\rho\div\delta x$,
see section 2 and \cite{L1}. 
The initial data for the linearized equations are
$$
\delta x\big|_{t=0}=\delta f_0,\qquad D_t\, \delta x\big|_{t=0}=\delta V_0.\tag 1.18
$$
The initial data are as before subject to constraints.
Let $\delta \tilde{x}$ be the formal power series solution in time $t$ to
(1.17)-(1.18). The compatibility condition of order $m$ is
$$
D_t^j\div\delta \tilde{x}=0,\qquad j=0,...,m-1.\tag 1.19
$$
The main difference between (1.17) and (1.14) is the higher order term
$\pa_i\big((\pa_k h)\delta x^k\big)$, since the term $\pa_i \delta h$,
depending on $\div \delta x$,  in (1.17) corresponds to the term $\pa_i h$,
depending on $\det{(\pa x/\pa y)}$ in (1.14). 
If we take $\overline{x}$ above to be a family of solutions 
of (1.14) depending on the parameter $r$, then our estimates below
show that a small change of initial conditions only give rise to a small 
change of the solution in Sobolev spaces. 
Our main result is
the following linearized stability result:
\proclaim{Theorem {1.}1} Let $\Omega$ be the unit ball in $\bold{R}^n$
and suppose that $x$ is a smooth solution of (1.14)
satisfying (1.9) for $0\!\leq \!t\!\leq\! T\!$.
Suppose that $(\delta f_0,\delta V_0)\!$ are smooth
satisfying the compatibility conditions of all orders $m$, 
i.e. (1.19) holds for all $m$. 
Then the linearized equations (1.17) have a smooth solution $\delta x$
for $0\leq t\leq T$ satisfying the initial conditions  (1.18).
Let ${\Cal N}$ be the exterior unit normal to
$\pa{\Cal D}_t$ parametrized by $x(t,y)$ and let
$\delta x_{\Cal N}={\Cal N}\cdot\delta x$ be the normal component of $\delta x$. Set
$$
E_r(t)=\| D_{t\,} \delta x(t,\cdot)\|_{H^r(\Omega)} +\|\delta x(t,\cdot)\|_{H^r(\Omega)}
+\|\div \delta x(t,\cdot)\|_{H^r(\Omega)}
+\|\delta x_{\Cal N}(t,\cdot)\|_{H^r(\pa \Omega)} \tag 1.20
$$
where $H^r(\Omega)$ and $H^r(\pa\Omega)$ are the Sobolev spaces in $\Omega$
respectively on $\pa\Omega$. Then there are constants $C_r$ depending only on $x$,
$r$ and $T$ such that
$$
E_r(t)\leq C_r E_r(0),\qquad\text{for}\qquad 0\leq t\leq T ,\qquad r\geq 0. \tag 1.21
$$

Furthermore, let
$\tilde{N}^r(\Omega)$ be the completion of $C^\infty(\overline{\Omega})$
in the norm $\|\delta x\|_{H^r(\Omega)}
+\|\div \delta x(t,\cdot)\|_{H^r(\Omega)}+\|\delta x_{\Cal N}\|_{H^r(\pa \Omega)} $.
 Then if the constraints in (1.19) hold for all orders $m$ and
$$
(\delta f_0,\delta V_0)\in
\tilde{N}^r(\Omega)\times H^r(\Omega)\tag 1.22
$$
 it follows that (1.17)-(1.18) has a solution
$$
(\delta x, D_t\delta x)\in C([0,T],\tilde{N}^r(\Omega)\times H^r(\Omega)).\tag 1.23
$$
\endproclaim
As we have argued, any smooth solution of (1.1)-(1.7) 
with ${\Cal D}_0$ diffeomorphic to the unit ball can be reduced to a 
smooth solution of (1.14) where $\Omega$ is the unit ball.
That there are initial data (1.18) such that (1.19)
hold for all $m$
follows by taking $\delta f_0$ and $\delta V_0$ compactly supported in the interior
of $\Omega$. The term $\|\delta x_{\Cal N}\|_{H^r(\pa\Omega)}$
is equivalent to the variation of the second fundamental form 
$\theta=\overline{\partial} {\Cal N}$ of the free boundary $\pa{\Cal D}_t$
measured in $H^{r-2}(\pa\Omega)$. For a general component we can only say that 
$\delta x\!\in\!\! H^{r-1/2}(\pa\Omega)$. The energy estimate (1.21) also hold in the
incompressible case when $\div\delta x\!=\!0$, see \cite{L1}, and in \cite{CL} we obtained 
similar bounds for $\|v\|_{H^r(\Omega)}+\|\theta\|_{H^{r-2}(\pa\Omega)}$
in the nonlinear incompressible case. 

Let us now outline the main ideas in the proof.
We will rewrite the linearized equations (1.17) in a geometrically invariant way
and use this to obtain energy bounds and existence.
We have defined our vector fields as functions of the Lagrangian
coordinates  $(t,y)\in[0,T]\times\Omega$ but we can also think of them
as functions of the Eulerian coordinates
$(t,x)\in {\Cal D}$, and we will make this identification
without saying that we compose with the inverse of the change of coordinate
$y\to x(t,y)$. The time derivative has a simple expression in the
Lagrangian coordinates but the space derivatives have a simpler expression
in the Eulerian coordinates, see (1.11). For the most part we will think of our
functions and vector fields in the Lagrangian frame but we use the
inner product coming from the Eulerian frame, i.e. in the Lagrangian frame
we use the pull-back metric of the Euclidean inner product:
$$
X\cdot Z=\delta_{ij} X^i Z^j=g_{ab} X^a Z^b,\qquad\text{where}\quad
X^a=X^i\frac{\pa y^a}{\pa x^i},\qquad
 g_{ab}=\delta_{ij}\frac{\pa x^i}{\pa y^a}\frac{\pa x^j}{\pa y^b}.\tag 1.24
$$
Here $X^i$ refers to the components of the vector $X$ in the Eulerian frame,
$X^a$ refers to the components in the Lagrangian frame, $g_{ab}$ is the metric
in the Lagrangian frame and $\delta_{ij}$ is the Euclidean metric. 
The letters $a,b,c,d,e,f,g$ will refer to indices in the Lagrangian frame whereas
$i,j,k,l,m,n$ will refer to the Eulerian frame.
The norms and most of the operators we consider have an invariant interpretation
so it does not matter in which frame they are expressed.
In the introduction we express the vector fields in the Eulerian frame but
later we express them in the Lagrangian frame.
The $L^2$ inner product is
$$
\langle X,Z\rangle=\int_{\dt} X\cdot Z \, dx=
\int_{\Omega} X\cdot Z\,\kappa\, dy=
\int_{\Omega} X\cdot Z\,\,\rho^{-1} k(y)\, dy,
\tag 1.25
$$
where $\kappa=\det{(\pa x/\pa y)}$ and $k=\kappa\rho\big|_{t=0}$.

Let us first point out that the boundary condition $\rho\big|_{\pa\Omega}=\overline{\rho}_0$
leads to that the energy is conserved for a solution of Euler's equations (1.12).
Let $Q(\rho)=\int_{\overline{\rho}_0}^\rho 2q(\rho)\rho^{-2}\, d\rho$,
where $q(\rho)=p(\rho)-p(\overline{\rho}_0)$.
Since $D_t(\rho\kappa)=0$ and $\rho \,D_t v_i=-\pa_i p(\rho)=-\pa_i q(\rho)$
it follows from the (1.25) and the divergence theorem that
$$\multline
\frac{d}{dt} \int_\dt\!\! (|V|^2  + Q(\rho))\rho \, dx=\!\int_\dt\!\!\!  D_t( |V|^2 +Q(\rho))
\rho \, dx
=\!\int_\dt\!\!\!\big( -2 V^i\pa_{i\,} q(\rho)+2q(\rho)\rho^{-1} D_t\rho\big)\, dx\\
=2\!\int_\dt \!\!\big(\div V\, q(\rho)+ q(\rho)\rho^{-1}D_t\rho\big)\, dx
-2\!\int_{\pdt} \!\!\!\! V_{\Cal N\,} q(\rho)\, dS=0,\!\!\!\!\!\endmultline
\tag 1.26
$$
where $V_{\Cal N}\!={\Cal N}_i V^i\!$ is the normal component and we also used (1.12).

We will obtain similar energy estimates for the linearized equations (1.17)
for energies containing an additional boundary term.
We will first rewrite the linearized equations  in a geometrically
invariant way. The last term in the first equation in (1.17) is
a positive symmetric operator  in the energy inner product
on vector fields satisfying the boundary condition
$\div X\big|_{\pa{\Cal D}_t}=0$:
$$
{C} X
=-\nabla \big( h^\prime(\rho) (\,\rho\, \div X+(\pa_k \rho) X^k)\big)
=- \nabla \big( h^\prime(\rho) \div \,(\rho X)\big),
\qquad\text{where} \quad  \nabla^i=\delta^{ij}\pa_j
\tag 1.27
$$
 i.e.
$\langle X,\rho\,{C}X\rangle \geq 0$ and
$\langle Z, \rho\,{C}X\rangle=\langle \rho\, {C}Z, X\rangle$
if $\div X\big|_{\pa{\Cal D}_t}=\div Z\big|_{\pa{\Cal D}_t}=0$. In fact,
if $\div X\big|_{\pa{\Cal D}_t}=0$ 
then $ h^\prime(\rho)\div\, (\rho X)\big|_{\pa{\Cal D}_t}=X^k\pa_k h\big|_{\pa{\Cal D}_t}
=  X_N\rho^{-1} \na_N p\big|_{\pa{\Cal D}_t}$ 
and integrating by parts we get 
$$
\langle Z, \rho \,{C}X\rangle =\int_{{\Cal D}_t}
\div\,(\,\rho  Z)\, \div\,(\rho X)\,
h^\prime(\rho)\, dx +\int_{\pa{\Cal D}_t} Z_N X_N (-\na_N p)\, dS\tag 1.28
$$
which proves the symmetry and the positivity follows
since $h^\prime(\rho)\geq c_1>0$ and $\na_N p\leq -c_0<0$.

We will also replace the time derivative by a time derivative that
preserves the boundary condition. Let 
$$
{\Cal L}_{D_t} X^i=D_t X^i-(\pa_k V^i)X^k
=\frac{\pa x^i}{\pa y^a}D_t  \Big(\frac{\pa y^a}{\pa x^k} X^k\Big),
\tag 1.29
$$
be the space time Lie derivative with respect to $D_t=(1,V)$
restricted to the space components or equivalently, the time derivative
of the vector field expressed in the Lagrangian frame. Let
$$
\hat{\Cal L}_{D_t} X^i={\Cal L}_{D_t} X^i+\div V\,  X^i=
\kappa^{-1} \frac{\pa x^i}{\pa y^a}D_t  \Big(\kappa\frac{\pa y^a}{\pa x^k} X^k\Big)
\tag 1.30
$$
be the modified Lie derivative that preserves the boundary condition,
$\div X\big|_{\pa\Omega}=0$. In fact
$$
\div\hat{\Cal L}_{D_t} X=\hat{D}_t \div X,
\qquad\quad\text{where}\qquad \hat{D}_t=D_t+\div V.\tag 1.31
$$

The linearized equations (1.17) can now be written as an evolution equation for
the operator ${C}$:
$$
 \ddot{X}+ {C} X=
 B\big(X,\dot{X}\big),\qquad\qquad \div X\big|_{t=0}=0,
\tag 1.32
$$
where $ X=\delta x$, $\dot{X}=
\hat{\Cal L}_{D_t} X$,
$\ddot{X}=\hat{\Cal L}_{D_t}^2 X$
and $B$ is a linear form with coefficients depending on $x$ and $\rho$.
 Associated with (1.32) is the energy
$$
E(t)=\langle \dot{X}, \rho\,\dot{X}\rangle
+\langle X,\rho \,({C}+I) X\rangle,\tag 1.33
$$
and we prove that 
$E^\prime\! \leq \! CE$ which gives the bound (1.21) for $r\!\!=0$.
The boundary term in (1.20) comes from (1.28). 
To obtain estimates for higher order derivatives one can apply
modified Lie derivatives with respect to tangential vector fields
as in \cite{L1}. This does however not prove existence for (1.32)
which is non-standard since ${C}$ is non-elliptic, time dependent
and the boundary condition is non-trivial. 

We use the orthogonal
projection onto divergence free vector fields  in the inner product
(1.25) to obtain
an equation for the divergence and an equation for the divergence free part.
The equations decouples
to highest order, and existence and estimates for the system follows from
existence and estimates for each equation with an inhomogeneous term.
The orthogonal projection is 
$$
 PX=X -\nabla q ,\qquad \text{where}
\qquad \triangle q=\div X,\qquad
q\Big|_{\Omega}=0. \tag 1.34
$$
We will obtain a system of equations for  $X_0=PX$ and $X_1=(I-P)X$
by projecting the linearized equation (1.32) onto divergence free vector
fields respectively the orthogonal complement.
Taking the divergence of (1.32) gives a wave equation for $\div X$
with Dirichlet boundary condition:
$$
\hat{D}_t^2 \div X_1-\triangle \big( h^\prime(\rho)\rho \div X_1\big)
=\triangle\big( X^k\pa_k h\big)+\div B(X,\dot{X}),\qquad
\quad \div X_1\big|_{\pa\Omega}=0
\tag 1.35
$$
for which existence is known if $h^\prime(\rho)\rho$
and the metric $g_{ab}$, hidden in 
$\triangle=\sum_{i=1}^n\pa_i^2=\kappa^{-1}\pa_a \,\kappa g^{ab}\pa_b$, 
are bounded from above and below and the right
hand side is thought of as a known function, see section 6. 
$X_1$ is then determined from $\div X_1$ by solving the Dirichlet problem:
$$
X_1=\nabla q_1,\qquad\quad \triangle q_1=\div X_1
,\qquad q_1\Big|_{\Omega}=0. \tag 1.36
$$

To obtain an equation for the divergence free part $X_0$
we project (1.32) onto divergence free vector fields.
It follows from (1.27) that 
$$
AX=P\,{C}X=P\big(-\nabla (X^k \pa_k h)\big),\tag 1.37 
$$
since $\div X\big|_{\pa\Omega}=0$ 
and the projection of a gradient of a function
that vanishes on the boundary vanishes
The operator $A$ is a positive symmetric operator on divergence free vector
fields, if condition (1.9) holds:
$$
\langle X,AZ\rangle=-\int_{\dt} \!X^i\pa_i (Z^k \pa_k h)\, dx=
\int_{\pdt}\!\! X_N Z_N \, (-\na_N h) \, dS,
\qquad\text{if}\qquad \div X=\div Z=0.\tag 1.38
$$
Furthermore, we note that the commutator of time derivatives with the projection is lower order
$$
[\hat{\Cal L}_{D_t},P]X^i=-P\big((\hat{\Cal L}_{D_t} \delta^{ij})\delta_{jk}(I-P)X^k\big)
\tag 1.39
$$
which follows since  $\hat{\Cal L}_{D_t}$ preserves the divergence free condition
and the projection of $\delta^{ij}\pa_j D_t q$ vanishes if $q\big|_{\pa\Omega}=0$
and hence $D_t\, q\big|_{\pa\Omega}=0$. Hence $P\ddot{X}_1=PB_2(X_1,\dot{X}_1)$
can be determined in terms of ${X}_1$ and $\dot{X}_1$.
Projection of (1.32) therefore gives an evolution equation for
the operator $A$ for the divergence free part:
$$
\ddot{X}_0+A X_0=-PB_2(X_1,\dot{X}_1) -AX_1+P B(X,\dot{X})
 \tag 1.40
$$
Existence for (1.40) with the right hand side thought of as a known function
is a generalization of the existence proof in \cite{L1}.
For the divergence part we have an equation which is equivalent to (1.35)-(1.36):
$$
\ddot{X}_1-\nabla \big( h^\prime(\rho)\rho\div X_1\big)
-PB_2(X_1,\dot{X}_1)
=(I-P)\nabla \big( (\pa_k h)X^k\big)+(I-P)B(X,\dot{X}),
\qquad\div X_1\big|_{\pa\Omega}\!\!=0
\tag 1.41
$$
We will show existence for the system (1.40)-(1.41)
for $(X_0,X_1)$
and from that we obtain a solution $X\!=\!X_0\!+\!X_1$ to (1.32), since these equations
are exactly the projection of (1.32) onto the divergence free vector fields
respectively the orthogonal component.
The system (1.40)-(1.41) can be solved by iteration.
If $X$ is an iterate, then from (1.35)-(1.36) and (1.40) 
we get $X_1$ and $X_0$ and a new iterate is $X_0\!+\!X_1$.
There is no loss of regularity in this procedure
since $\div X$ has the same space regularity as $X$.

However, in order for it to be possible to solve (1.35) the initial conditions
and the equation must be compatible with a formal power series solution
satisfying the boundary conditions, and we must make sure that this is true
at each step of the iteration.  (1.32) gives $\ddot{X}$ in terms
of only space derivatives of $X$ and $\dot{X}$ and we hence obtain a formal power
series solution in time, the first two terms coming from the initial conditions.
Our assumption is that this formal power series solution satisfies the
boundary condition. From this formal power series one can construct an approximate
solution $\tilde{X}$ satisfying the initial conditions, the equation to all orders
as $t\to 0$, and the boundary condition. We can then take the approximate solution
as our first iterate or equivalently subtract off the approximate solution from
$X$, which produces an inhomogeneous term vanishing to all orders as $t\to 0$
and vanishing initial conditions.

Let, us now conclude the introduction by giving the main estimates we use.
Since the time derivative preserves the boundary condition it is natural to
use norms which also contain time derivatives up to full order in the proof,
and the estimate in Theorem {1.}1 afterwards follow from these.
Let
$$
\|X(t)\|_{H^r}=\|X(t,\cdot)\|_{H^r(\Omega)},\quad
 {\|} X(t){\|}_r=\!\!\sum_{s+k\leq r}\!\| D_t^k X(t)\|_{H^s},
\quad \langle X_0(t)\rangle_r=\|X_{0N}(t,\cdot)\|_{H^r(\pa\Omega)}\tag 1.42
$$
where $X_{0N}=X_0\cdot N$ is the normal component.
For the divergence free equation:
$$
\ddot{X}_0+A X_0=F_0,\qquad PF_0=F_0\tag 1.43
$$
we have the estimate
$$
{\|}\dot{X}_0(t){\|}_r+{\|} X_0(t){\|}_r+\langle X_0(t)\rangle_r
\leq C\Big({\|}\dot{X}_0(0){\|}_r+{\|} X_0(0){\|}_r+\langle X_0(0)\rangle_r
 + \int_0^t {\|} F_0{\|}_r\, d\tau\Big).\tag 1.44
$$
This is a generalization of the estimate for the incompressible case in \cite{L1}.
For $r=0$, one uses the symmetry and positivity  (1.38) of $A$
to prove that $E=\langle \dot{X}_0,\dot{X}_0\rangle
+\langle X_0,(A+I)X_0\rangle $ satisfies $E^{\,\prime}\leq C E$. Note that 
the boundary term comes from using (1.38). For $r>0$,
it follows from commuting modified Lie derivatives with respect to
tangential vector fields through the equation to obtain similar equations and
estimates for these, together with better estimates
for the curl since the curl of $A$ vanishes.
For
$$
\ddot{X}_1-\nabla \big( h^\prime(\rho)\rho\div X_1\big)
-PB_2(X_1,\dot{X}_1)
=F_1,\qquad (I-P)F_1=F_1,\qquad\quad \div X_1\big|_{\pa\Omega}=0,\tag 1.45
$$
we have the estimate
$$
{\|} X_1(t){\|}_{r+1}\leq
C\Big({\|} X_1(0){\|}_{r+1}
+\int_0^t {\|}\dot{F}_1{\|}_{r-1}\, d\tau\Big).\tag 1.46
$$
The last estimate follows from estimating the wave equation (1.35) with
the right hand side replaced by $\div F_1$ and
inverting the Laplacian (1.36).
In (1.46) we do not need space derivatives up to highest
order of $F_1$, since one obtains space derivatives
from time derivatives through inverting the Laplacian in the
wave equation. Using that the right hand side of (1.41) is a gradient
(1.46) also holds for $r\!=\!0$.

With $F_0$ equal to the right hand side of (1.40)
and $F_1$ equal to the right hand side of (1.41) the norms in the
integrals in (1.44) and (1.46) can be estimated by the sum of the
norms in the left of (1.44) and (1.46) and this gives {\it{ a priori}}
bounds as well as uniform estimates for iterates, if we let the $F_0$ and $F_1$
be obtained from the previous iterate and solve (1.43) and (1.45) for the new
iterate. Once we obtained the solution to the system (1.40)-(1.41), $X=X_0+X_1$ is
the solution to (1.32). The norm in (1.20) is bounded by the sum of the norms
in the left of (1.44) and (1.46). Using (1.32) one can bound time derivatives
in terms of space derivatives and using that the projection is continuous in
the norms (1.42) it follows that the norms in the left of (1.44) and (1.46)
can be bounded by (1.20).

\head 2. Lagrangian coordinates and the linearized equation.\endhead
Let us introduce Lagrangian coordinates in which the boundary becomes fixed.
Let $\Omega$ be a domain in $\bold{R}^n$ and let
$f_0:\Omega\to {\Cal D}_0$ be a diffeomorphism.
We assumed that ${\Cal D}_0$ is diffeomorphic to the unit ball and 
that $v(t,x)$, $p(t,x)$, $(t,x)\in  {\Cal D}$ are given satisfying the boundary
conditions (1.3)-(1.4).
The Lagrangian coordinates $y$ are given by solving for the Eulerian coordinates
$x=x(t,y)=f_t(y)$ in 
$$
{d x}/{dt}=V(t,x(t,y)), \qquad x(0,y)=f_0(y),\quad y \in \Omega\tag 2.1
$$
Then $f_t:\Omega\to {\Cal D}_t$ is a diffeomorphism,
and the boundary becomes fixed in the new $y$ coordinates.
Let us introduce the notation
  $$\align
\vspace{-0.06in}
\quad D_t&=\frac{\partial }{\partial t}\Big|_{y=constant}=
\frac{\partial }{\partial t}\Big|_{x=constant}
+\,  V^k\frac{\partial}{\partial x^k},\tag 2.2\\
\vspace{-0.3in}
&
\endalign
$$
for the material derivative and
$$\align
\vspace{-0.07in}
\partial_i&=\frac{\partial}{\partial {x^i}}=\frac{\partial y^a}{\partial x^i}
\frac{\partial}{\partial y^a}. \tag 2.3
\endalign
$$
for the partial derivatives.

In  these coordinates Euler's equation (1.1) become
$$
\rho \, D_t^2 x_i+\pa_{i\,} p=0, \qquad (t,y)\in[0,T]\times\Omega\tag 2.4
$$
and the continuity equation (1.2) become
$$
D_t \rho+\rho\,\div V=0,\qquad (t,y)\in[0,T]\times\Omega\tag 2.5
$$
Here the pressure $p=p(\rho)$ is assumed to be a given smooth
strictly increasing function of the density $\rho$. Let
$\overline{\rho}_0$ be defined by  $p(\overline{\rho}_0)=0$. Let
$h$, the enthalpy, be defined by
$$
h(\rho)=\int_{\overline{\rho}_0}^\rho
p^\prime(\rho)\rho^{-1}\, d\rho.\tag  2.6
$$
Then (2.4) becomes
$$
D_t^2 x_i+\pa_i h=0,\quad (t,y)\in[0,T]\times\Omega.\tag 2.7
$$

The density $\rho$ satisfies (2.5)
but since
$$
\kappa=\det{(\pa x/\pa y)},\tag 2.8
$$
satisfies
$$
D_t\kappa-\kappa \div V=0\tag 2.9
$$
 it follows
that $\rho\!=\!\rho_0 \kappa_0/\kappa$, where $\rho_0$ and $\kappa_0$ are the initial
values. By a theorem of \cite{DM} one can arbitrarily prescribe the
volume form $\kappa_0$ up to a constant so we take $\kappa_0\!=\!k/\rho_0$,
where $k$ is a constant,
and $\Omega$ to be the unit ball, by composing with a diffeomorphism, 
since we assumed that ${\Cal D}_0$ is diffeomorphic to a unit ball. 
 Hence $\rho$ is determined from $x$:
$$
\rho=k/\kappa\tag 2.10
$$
By choosing the constant $k$ appropriately the boundary condition
(1.4) can hence be expressed
$$
\kappa\big|_{\pa\Omega}=1\tag 2.11
$$

Since $h$ is a function of $\rho$ which in turn by our choice
(2.10) is a function of
$\kappa=\det{(\pa x/\pa y)}$ we can think of $h$ as a function of $\kappa$.
(2.7) is then an equation involving the coordinate $x$ only and initial data
for $\rho$ is included in the choice of initial mapping $f_0$.
Initial data for (2.7) are
$$
x\big|_{t=0}=f_0,\qquad D_t\, x\big|_{t=0}=V_0\tag 2.12
$$
In order for (2.7)
to have a smooth solution satisfying (2.11), initial data has to satisfy the constraints
$\det{(\pa f_0/\pa y)}\big|_{\pa\Omega}\!\!=1$ and
$\div V_0\big|_{\pa\Omega}\!\!=0$, by (2.9).
Taking the divergence of (2.7) gives
$$
D_t\div V+(\pa_i V^k)\pa_k V^i +\triangle h=0\tag 2.13
$$
which leads to further conditions.
Since (2.7) gives $D_t^2 x$ in terms of space derivatives of $x$
we can obtain a formal power series solution in time $t$, $\tilde{x}$, to
(2.7) satisfying the initial conditions (2.12). The compatibility
condition of order $m$ is the requirement
that the formal power series solution up to terms of order $m$
satisfy the boundary condition in (2.11):
$$
D_t^j \big(\det{(\pa \tilde{x}/\pa y)}-1\big)\big|_{0\times\pa\Omega}=0,
\qquad j=0,...,m-1\tag 2.14
$$

At this point we also remark that we get a wave equation for $h$.
Since $h$ is a strictly increasing function of $\rho$ we can think of
$\rho=\rho(h)$ as function of $h$.
Hence with $e(h)=\ln{\rho(h)}$ (2.5) instead become
$$
D_t\, e(h)+\div V=0. \tag 2.15
$$
and this together with (2.13) gives a wave equation for $h$
with Dirichlet boundary conditions:
$$
D_t^2 e(h)-\triangle h-(\pa_i V^k)\pa_k V^i=0,\qquad
h\Big|_{\pa \Omega}\!\!\!=0.\tag 2.16
$$
Here
$$
\triangle h\!=\!\sum_{i}\pa_i^2 h\!=\!
\kappa^{-1}\pa_a\big(\kappa g^{ab}\pa_b h\big),\qquad\text{where}\quad
g_{ab}=\delta_{ij}\frac{\pa x^i}{\pa y^a}\frac{\pa x^j}{\pa y^b},\tag 2.17
$$
is the metric in the Lagrangian coordinates and $g^{ab}$ is its inverse.
Here $\pa_a =\pa/\pa y^a $ and we use
the convention that differentiation with respect to the Eulerian
coordinates is denotes by letters $i,j,k,l,m,n$ and with respect to the
Lagrangian coordinates is denoted by $a,b,c,d,e,f$.
In order for (2.16) to be solvable we must have that
$$
0< e^\prime+1/e^\prime \leq c_1
\qquad \sum_{a,b=1}^{n} \big(| g^{ab}|+|g_{ab}|\big)\leq n \,c_1^2,\qquad
 |\pa x/\pa y|^2+|\pa y/\pa x|^2\leq c_1^2\tag 2.18
$$
for some constant $0<c_1<\infty$. The first condition is related to that 
the pressure
is assumed to be a strictly increasing smooth function of the density.
The second and third condition are equivalent and says that
the coordinate mapping is a diffeomorphism.
Furthermore, it is well-known that one needs compatibility conditions to solve (2.16).

\comment
and the boundary condition (1.4) become
$$
h\big|_{\pa\Omega}=0.\tag 2.7
$$
Since $\pa_i h=\rho^{-1}\pa_i p$, Christodoulou's physical condition (1.9) become
$$
\na_{\Cal N}\, h\leq -c_0 <0,\quad\text{on}\quad \partial {\Omega},
\qquad \text{where}\quad \na_{\Cal N}={\Cal N}^i\partial_{x^i}. \tag 2.8
$$
Since $h$ is a strictly increasing function of $\rho$ we can think of
$\rho=\rho(h)$ as function of $h$.
Hence with $e(h)=\ln{\rho(h)}$ (2.5) instead become
$$
D_t\, e(h)+\div V=0. \tag 2.10
$$
Note that this, together with (2.7) also implies that $\div V\big|_{\pa\Omega}=0$.
We have now reduced the system (1.1)-(1.2) to (2.6) and (2.10)
and we have reduced the boundary condition (1.3) to (2.7).
We have reduced the physical condition (1.6) to (2.8).
However, as we will see in below, in order for it to be
solvable there are also a couple of other conditions that enter
(2.13)-(2.14).

Taking the divergence of (2.6) using (2.10), and that
$[D_t,\pa_i]=-(\pa_i V^k)\pa_k$, gives a wave equation for $h$
with Dirichlet boundary conditions:
$$
D_t^2 e(h)-\triangle h-(\pa_i V^k)\pa_k V^i=0,\qquad
h\Big|_{\pa \Omega}\!\!\!=0.\tag 2.11
$$
Here
$$
\triangle h\!=\!\sum_{i}\pa_i^2 h\!=\!
\kappa^{-1}\pa_a\big(\kappa g^{ab}\pa_b h\big),\qquad\text{where}\quad
g_{ab}=\delta_{ij}\frac{\pa x^i}{\pa y^a}\frac{\pa x^j}{\pa y^b},\tag 2.12
$$
is the metric in the Lagrangian coordinates and $g^{ab}$ is its inverse.
Here $\pa_a =\pa/\pa y^a $ and we use
the convention that differentiation with respect to the Eulerian
coordinates is denotes by letters $i,j,k,l,m,n$ and with respect to the
Lagrangian coordinates is denoted by $a,b,c,d,e,f$.

In order for (2.11) to be solvable we must have that
$$
0< e^\prime+1/e^\prime \leq c_1
\qquad \sum_{a,b=1}^{n} \big(| g^{ab}|+|g_{ab}|\big)\leq n \,c_1^2,\qquad
 |\pa x/\pa y|^2+|\pa y/\pa x|^2\leq c_1^2\tag 2.13
$$
for some constant $0<c_1<\infty$. The first condition is related to that the pressure
is assumed to be a strictly increasing smooth function of the density.
The second and third condition are equivalent and it is a condition on the coordinates.

Without getting further into the physical motivation of the condition (2.13)
we remark that from a mathematical point of view it is necessary for the
well-posedness of the problem since one can not solve
the wave equation (2.11) without this condition.

In order for (2.11) to have a formal power series solution in
$C^{m-1}$ we must also assume that $m$ compatibility conditions on
initial data are satisfied. These are obtained as follows. Let
$(\tilde{x},\tilde{h})$ be a formal power series solution
$$
D_t^j h\big|_{\pa\Omega}=0, \qquad \text{when}\quad t=0,\quad \text{for}
\quad j=0,...,m-1.\tag 2.14
$$
For $j=0,1$ this simply means that the initial conditions should vanish
on the boundary, for $j=2$ this means that the sum of the last two terms
(2.11) should vanish on the boundary and for larger $j$ one gets more
conditions from taking time derivatives of (2.11).

Using that also $[\delta-\delta x^k\pa_k,\triangle]=0$ we get from applying
$\delta-\delta x^k\pa_k$ to (2.11)
$$
D_t^2 (e^\prime(h)\delta h)-\delta x^k\pa_k D_t^2 e(h)
-\triangle\big(\delta h-\delta x^k\pa_k h\big)
-2(\pa_i V^k)\pa_k(\delta V^i\!-\delta x^l\pa_l  V^i)=0 ,\qquad
\delta h\Big|_{\pa \Omega}\!\!\!=0.\tag 2.16
$$
where $\delta V^i=D_t\delta x^i$.
Since we are considering
a variation within families of solutions satisfying the $m$ compatibility
conditions (2.14) this then translates into the same $m$
compatibility conditions for the variation in order for
the last part of (2.16) to hold. These are
$$
D_t^j \delta h\big|_{\pa\Omega}=0, \qquad \text{when}\quad t=0,\quad \text{for}
\quad j=0,...,m-1.\tag 2.17
$$

Since
$$\multline
[\pa_i ,D_t^2]\delta x^i=D_t\big([\pa_i,D_t]\delta x^i\big)
+[\pa_i,D_t]D_t\delta x^i=
2(\pa_i V^k)\big(\pa_k\delta V^i-(\pa_k V^l)\pa_l \delta x^i\big)
+(\pa_i D_t V^k)\pa_k \delta x^i\\
=2(\pa_i V^k)\pa_k\big(\delta V^i-\delta x^l\pa_l V^i\big)
+\delta x^l\pa_l\big((\pa_i V^k)\pa_k V^i \big)
+\pa_k\big((\pa_i D_t V^k)\delta x^i\big)-(\pa_i \pa_k D_t V^k)\delta x^i\\
=2(\pa_i V^k)\pa_k\big(\delta V^i-\delta x^l\pa_l V^i\big)
+\pa_k\big((\pa_i D_t V^k)\delta x^i\big)
-(\pa_i D_t \div V)\delta x^i
\endmultline\tag 2.18
$$
Hence if we add (2.19) and (2.16) and use (2.10) we get
$$
D_t^2\big(\div \delta x+e^\prime(h)\delta h\big)=0,
\tag 2.20
$$
It follows that
$$
\div\delta x=-e^\prime\delta h,\tag 2.21
$$
Actually, this also follows from our original assumption that $h$ was
a function of the Jacobian $\kappa$ since if $\sigma=\ln \kappa$ then
$\delta \sigma=\div \delta x$, see \cite{L1}, so
$\delta h(\sigma)=h^\prime \delta \sigma=h^\prime \div\delta x=...=
-\div\delta x/e^\prime $.
\endcomment

Let us now derive the linearized equations.
The calculations that follows below are similar to those in \cite{L1}
since the equation (2.7) mathematically is the same as the equation for the
incompressible case with the enthalpy $h$ replaced by the pressure $p$.
We therefore refer the reader to \cite{L1} for more details.

We now assume that we have a smooth solution $x=x(t,y)$ of (2.7)
satisfying (1.9) for $0\leq t\leq T$
and we will derive the linearized equations at this solution.
Assume that $\overline{x}\!=\overline{x}(t,y,r)$ is a smooth function also
of the extra parameter $r$ such that $\overline{x}\big|_{r=0}\!\!=x$
and set $\delta x=\pa \overline{x}/\pa r\big|_{r=0}$.
Then the linearized equations are the requirements on $\delta x$ that
$\overline{x}$ satisfies (2.7) and (2.10)-(2.11)
up to terms bounded by $r^2$ as $r\!\to\! 0$.
Let $\delta$ be a variation
in the Lagrangian coordinates, i.e. a derivative
$$
\delta f=\pa f/\pa r\big|_{r=0}.\tag 2.19
$$
 with respect to the
parameter $r$ when $t$ and $y$ are fixed.
Then $[\delta,D_t]=0$,
$$
[\delta,\pa_i]\!=\!-(\pa_i \delta x^k)\pa_k,\tag 2.20
$$
so
$[\delta\! -\delta x^k\pa_k,\pa_i]\!=\!0$. Applying $\delta-\delta x^k\pa_k$
to (2.7) gives:
$$
D_t^2 \delta x_i-(\pa_k D_t^2 x_i)\delta x^k
- \pa_i \big( \delta x^k\pa_k h- \delta h \big)=0\tag 2.21
$$
Since $h=h(\rho)$ where $\rho=k/\kappa$ and
$$
\delta \kappa=\kappa \div\delta x. \tag 2.22
$$
it follows that
$$
\delta h=-h^\prime(\rho)\rho\div \delta x.\tag 2.23
$$
The variation of the boundary condition (2.11) become
$$
\div\delta x\big|_{\pa\Omega}=0.\tag 2.24
$$
The initial data for (2.21) with $\delta h$ given by (2.23) are
$$
\delta x=\delta f_0,\qquad\quad D_t \delta x=\delta V_0\tag 2.25
$$
In order for it to be possible to have a smooth solution of (2.21)
and (2.23)-(2.24) initial data (2.25) must satisfy certain compatibility conditions.
The initial data are subject to the constraints
$\div \delta f_0\big|_{\pa\Omega}=0$ and
$\div \delta V_0=(\pa_i V^k)\pa_k\delta x^k$.

Taking the divergence of (2.21) using (2.23) and (2.24) gives a wave
equation for $\div\delta x$ with Dirichlet boundary conditions:
$$
D_t^2 \div \delta x-\delta x^i\pa_i D_t\div V
-\triangle  \big(\delta x^k\pa_k h- \delta h \big)
+2(\pa_i V^k)\pa_k(\delta V^i\!-\delta x^l\pa_l V^i)=0,\tag 2.26
$$
which gives further  conditions.
Since (2.21) gives $D_t^2 \delta x$ in terms of only space derivatives
of $\delta x$, this gives a formal power series solution in time $t$, which we call
$\delta \tilde{x}$. The compatibility condition of order $m$ is the requirement that
the formal power series solution satisfies the boundary condition (2.24):
$$
D_t^j \, \div\delta\tilde{x}=0,\qquad j=0,...,m-1.\tag 2.27
$$
The basic assumption in solving the system (2.21)-(2.25) is that
one should assume that $\div \delta x$ has the same space
regularity as $\delta x$.

Let us now express also the vector field in the Lagrangian frame.
Let
$$
 W^a=\frac{\pa y^a}{\pa x^i }\delta x^i\tag 2.28
$$
Then,
$$
D_t\, \delta x^i=D_t \big( W^b\pa x^i/\pa y^b\big)
= (D_t W^b)\pa x^i/\pa y^b+ W^b\pa V^i/\pa y^b
= (D_t W^b)\pa x^i/\pa y^b+ \delta x^k\pa_k V^i\tag 2.29
$$
and multiplying with the inverse $\pa y^a/\pa x^i$ gives
$$
D_t W^a= \frac{\pa y^a}{\pa x^i } {\Cal L}_{D_t} \delta x^i,
\qquad\text{and}\qquad
\hat{D}_t W^a= \frac{\pa y^a}{\pa x^i } \hat{\Cal L}_{D_t} \delta x^i.\tag 2.30
$$
where the Lie derivative and modified Lie derivative are given by (1.29)-(1.30)
and
$$
\hat{D}_t W^a=D_t W^a+(\div V)W^a=\kappa^{-1} D_t(\kappa W^a).\tag 2.31
$$
Since the divergence is invariant
$$
\div \delta x=\div W=\kappa^{-1}\pa_a\big(\kappa W^a\big)\tag 2.32
$$
it therefore follows that
$$
\div \hat{D}_t W=\hat{D}_t \div W.\tag 2.33
$$

Differentiating (2.30) once more gives
$$
D_t^2 \delta x^i-(\pa_k D_t V^i)\delta x^k=(D_t^2 W^b)\pa x^i/\pa y^b
+2(D_t W^b)\pa V^i/\pa y^b\tag 2.34
$$
It follows that
$$\aligned
\frac{\pa x^i}{\pa y^a }\big(D_t^2 \delta x^i-(\pa_k D_t V^i)\delta x^k\big)
&=\frac{\pa x^i}{\pa y^a }\frac{\pa x^i}{\pa y^b }D_t^2 W^b
+2(D_t W^b)\frac{\pa x^i}{\pa y^b }\frac{\pa x^j}{\pa y^a }
\pa_i v_j\\
&=g_{ab} D_t^2 W^b+ (D_t\, g_{ab}-\omega_{ab}) D_t W^b
\endaligned\tag 2.35
$$
where $g_{ab}$ is given by (2.17) and
$$
D_t\, g_{ab}=\frac{\pa x^i}{\pa y^a}\frac{\pa x^j}{\pa y^b}\big(\pa_i v_j\!+\pa_j v_i\big),
\qquad \omega_{ab}
=\frac{\pa x^i}{\pa y^a}\frac{\pa x^j}{\pa y^b}\big(\pa_i v_j\!-\pa_j v_i\big).
\tag 2.36
$$
With $\pa_a=\pa/\pa y^a$ the linearized equation (2.21) and (2.23) become
$$
g_{ab} D_t^2 W^b-\pa_a\big( (\pa_c h) W^c
-\delta h\big)
=-\big(D_t{{g}}_{ac}-\omega_{ac} \big) D_t W^c, \qquad \delta h=-p^\prime\, \div W
\qquad \tag 2.37
$$

Let $\hat{D}_t$ be as in (2.31), i.e.
$\hat{D}_t=(D_t+\dot{\sigma})$,  where $\sigma=\ln{\kappa}$ and
$\dot{\sigma}=D_t\sigma=\div V$. Then
$$
D_t^2=\hat{D}_t^2-2\dot{\sigma}\hat{D}_t+\dot{\sigma}^2-\ddot{\sigma},
\qquad D_t=\hat{D}_t-\dot{\sigma},\qquad \ddot{\sigma}=D_t^2\sigma.\tag 2.38
$$
Hence, with $\dot{W}=\hat{D}_t W$ and $\ddot{W}=\hat{D}_t^2 W$,
we can write (2.37) as $LW=0$, where 
$$
 LW=\ddot{W}^a-g^{ab}\pa_b\big( (\pa_c h) W^c-\delta h\big)
 -B^a(W,\dot{W}),\qquad \delta h=-p^\prime \,\div W\,\tag 2.39
$$
where
$$
B^a(W,\dot{W})=-g^{ab}\big(\dot{{g}}_{bc}-\omega_{bc} \big)
(\dot{W}^c-\dot{\sigma} W^c)
+2\dot{\sigma} \dot{W}^a+(\ddot{\sigma}-\dot{\sigma}^2 )W^a.\tag 2.40
$$

\head 3. The compatibility conditions, statement of the theorem
and the lowest order energy estimate.\endhead
We now consider the linearized operator
$$
 LW=\ddot{W}+{C} W-B(W,\dot{W})\tag 3.1
$$
where $\dot{W}=\hat{D}_t W$,  $\ddot{W}=\hat{D}_t^2 W$, $\hat{D}_t=D_t+(\div V)$,
$B$ is the bounded operator
given by (2.40) and
$$
{C} W^a=-g^{ab}\pa_b\big( (\pa_c h) W^c+p^\prime\div W\big).\tag 3.2
$$
We want to show existence and estimates for the linearized
equations with an inhomogeneous term $F$,
$$
L W=F,\tag 3.3
$$
with initial data
$$
 W\big|_{t=0}=\tilde{W}_0,\quad \dot{W}\big|_{t=0}=\tilde{W}_1, \tag 3.4
$$
and boundary data
$$
\div W\big|_{\pa \Omega}=0.\tag 3.5
$$

The reason for the inhomogeneous term $F$ is that one can reduce to the case of vanishing
initial data and an inhomogeneous term $F$ that vanishes to all orders as $t\to 0$
and it is easier to first prove existence for this case.
Differentiating (3.3) with respect to time we get
$$
\hat{D}_t^{k+2} W=B_k\big( W,.,\hat{D}_t^{k+1}W,\pa W,...,
\pa \hat{D}_t^{k} W,
\pa^2 W,...,\pa^2 \hat{D}_t^k W\big)+\hat{D}_t^k F,\tag 3.6
$$
for some function $B_k$.
Let us therefore define functions of space only by
$$
\tilde{W}_{k+2}=B_k\big(\tilde{W}_0,...,\tilde{W}_{k+1},
\pa \tilde{W}_0,...,\pa  \tilde{W}_k,
\pa^2  \tilde{W}_0,...,\pa^2  \tilde{W}_k\big)\big|_{t=0}+\hat{D}_t^k F\big|_{t=0},
\qquad k\geq 0   \tag 3.7
$$
In view of (3.5) it follows that $0=\hat{D}_t^k \div W\big|_{\pa\Omega}
=\div\hat{D}_t^k W\big|_{\pa\Omega}$ so we must have
$$
\div  \tilde{W}_k\big|_{\pa\Omega}=0,\qquad\quad  k=0,...,m\tag 3.8
$$
(3.8) is called the $m^{th}$ order compatibility condition and in order for
it to be possible for (3.3)-(3.5) to have a smooth solution these have to hold for
all orders $m$.

We now define the approximate power series solution by
$$
\tilde{W}(t,y)=\frac{\kappa(0,y)}{\kappa(t,y)}
\sum_{k=0}^{\infty}\chi(t/\varepsilon_k)  \tilde{W}_k(y) t^k/k! .\tag 3.9
$$
Here $\chi$ is smooth $\chi(s)=1$ for $|s|\leq 1/2$ and $\chi(s)=0$ for $|s|\geq 1$.
The sequence $\varepsilon_k>0$ can be chosen so that the series
converges in $C^m([0,T],H^m)$ for any $m$
if  take $(\|\tilde{W}_k\|_{H^k}+1)\varepsilon_k\leq 1/2$.
It follows that
$$
\div \tilde{W}\big|_{\pa\Omega}=0.\tag 3.10
$$
(3.9) multiplied with $\kappa(t,y)$ is a power series expansion of $\kappa W$
and it hence follows that (3.9) satisfies (3.3)-(3.5) to all orders as $t\to 0$:
$$
D_t^k\big( L \tilde{W}-F\big)\big|_{t=0}=0,\qquad\quad k=0,...\tag 3.11
$$
It follows that we can reduce (3.3)-(3.5) to the case with vanishing initial data
and an inhomogeneous term that vanishes to all orders as $t\to 0$,
by replacing $W$ by $W-\tilde{W}$ and $F$ by $F-L \tilde{W}$ in (3.3).

Let us now introduce some notation:
\demo{Definition {3.}1} Let
$$
\|W(t)\|_{H^r}=\| W(t,_{\!}\cdot)\|_{H^r(\Omega)}.\tag 3.12
$$
and
$$
{\|}W(t){\|}_r=\!\!\sum_{s+k\leq r}\!\!\|\hat{D}_t^k W(t)\|_{H^s}.\tag 3.13
$$
Let $N$ be the exterior unit normal to $\pa\Omega$ in the metric $g_{ab}$,
or equivalently, $N^a={\Cal N}^i \pa y^a/\pa x^i$. Set
$$
\langle W(t)\rangle_{r}\!=\|W_N(t,_{\!}\cdot)\|_{H^r(\pa \Omega)}.
\tag 3.14
$$
where $W_N=W\cdot N$ is the normal component.
\enddemo

\proclaim{Theorem {3.}1} Suppose that $p=p(\rho)$ is a strictly increasing smooth function
of $\rho$. Suppose also that $x$ is a smooth solution of
(2.7), such that (1.9) hold for $0\leq t\leq T$. 
Suppose that the inhomogeneous term $F$ in (3.3) is smooth for
$0\leq t\leq T$. Suppose also that the initial conditions (3.4)
are smooth and satisfy the $m^{th}$ order compatibility conditions (3.8),
for all $m=0,1,...$.
 Then the linearized equations (3.3)-(3.5)
 have a smooth solution for $0\leq t\leq T$.

Let
$$
\tilde{E}_r(t)={\|}\dot{W}(t){\|}_{r}+{\|} W(t){\|}_{r}
+\langle W(t)\rangle_{r}
+{\|}\div W (t){\|}_{r},\tag 3.14
$$
where $\dot{W}=\hat{D}_t W=D_t W+(\div V) W$.
Then there is a constant $C$ depending only on $x$, $r$ and $T$
such that for $0\leq t\leq T$ we have
$$
\tilde{E}_r(t)\leq C\Big(\tilde{E}_r(0)+\int_0^t {\|} F{\|}_r\, d\tau \Big).\tag 3.15
$$
\endproclaim

Theorem {1.}1 follows from Theorem {3.}1 since the norm (3.14) is equivalent to
$$
{E}_r(t)=\|\dot{W}(t)\|_{H^r}+\| W(t)\|_{H^r}+\langle W(t)\rangle_{H^r}
+\|\div W (t)\|_{H^r},\tag 3.16
$$
if $F$ vanishes.
In fact, by (3.6) one can express time derivatives in terms of space derivatives
of the same order or less and using induction it follows that
$$
E_r\leq \tilde{E}_r\leq C_r\big( E_r +{\|}F{\|}_{r-1}\big)\tag 3.17
$$

In this section we show the lowest order energy estimates for an equation of the form
$$
\ddot{W}+{C}W=B(W,\dot{W})+F\tag 3.18
$$
where $\dot{W}=\hat{D}_t W=\kappa^{-1}D_t(\kappa W)$, $\ddot{W}=\hat{D}_t^2 W$,
$$
{C} W^a=-g^{ab}\pa_{\, b} \big( p^{\,\prime}\big(\div
W+(\pa_c e) W^c\big)\big)= -g^{ab}\pa_{\, b} \big( h^\prime \div\,
(\rho\, W)\big), \qquad e(\rho)=\ln\rho,\quad p^\prime(\rho)
=h^\prime(\rho)\rho
\tag 3.19
$$
and $B$ is any bounded linear operator.
The energy is:
$$
\multline
E=\langle \dot{W},\rho\dot{W}\rangle +\langle W,\rho\,({C}+I)W\rangle\\
=\int_{\Omega} g_{ab} \dot{W}^a \dot{W}^b + g_{ab}{W}^a {W}^b
+p^\prime\big(\div\,(\rho_{\,}W)/\rho\big)^2  \, \rho\kappa dy
+\int_{\pa\Omega} W_N^2 \, (-\na_N p)\, dS.
\endmultline \tag 3.20
$$
Now, for any symmetric operator $B$ we have
$$
\frac{d}{dt}\langle W, \rho BW\rangle =\frac{d}{dt}
\int_{\Omega} \kappa W^a \rho \underline{B}W_{\!a }\, dy
=2\langle \dot{W}, BW\rangle
+\langle W,\rho B^\prime W\rangle,
\tag 3.21
$$
where $ \dot{W}=\kappa^{-1}D_t(\kappa W)$ and
$\rho B^\prime$ is the time derivative of the operator $\rho B$
considered as an operator from the vector fields to the one forms:
$$
\rho B^\prime W^a =g^{ab}( D_t( \rho \underline{B}W_b) -\rho \underline{B} \dot{W}_b ),
\qquad \quad \underline{B}W_b=g_{bc} B W^c,
\tag 3.22
$$
Since $\langle W,\rho W\rangle =\langle W,\rho\,  GW\rangle$, where $G=I$,
it follows that
$$
\dot{E}=2\langle \dot{W},\rho\,\ddot{W}+\rho\,(C+I)W\rangle
+\langle \dot{W},\rho\, {G}^\prime\dot{W}\rangle
+\langle W,\rho\, ({C}^\prime+{G}^\prime)W\rangle.
\tag 3.23
$$
where $\rho\,  \underline{G}^\prime W_a=D_t\big(\rho\, g_{ab}\kappa) W^b $ and
$\underline{C}^\prime W_a= D_t\, {\underline{C}}W_a-{\underline{C}}\dot{W}_a
+\dot{e}\,{\underline{C}}W_a$,
where $\dot{e}=\dot{\rho}/\rho$.  Since
$D_t(\rho\,\kappa)\!=\!0$ and $D_t \big(\kappa \div \, (\rho W)\big)=
\kappa \div\big( D_t (\rho W)\big)$ we get 
$$\multline 
D_t\, \underline{C} W_a
=-\pa_a D_t\big( (\rho\kappa)^{-1} p^\prime(\rho) \kappa \div \,(\rho \,W)\big)\\
=-\pa_a\Big( p^{\prime\prime}(\rho)\, \dot{e} \,  \div \,(\rho W)
+ p^\prime(\rho) \rho^{-1} \div \,(\dot{\rho}\,  W+\rho\,\dot{W}) \Big)
\endmultline \tag 3.24 
$$
so
$$
 \underline{C}^\prime W_a=
-\pa_a\big( p^{\prime\prime}(\rho)\, \dot{e} \,  \div \,(\rho W)\big)
- \pa_a\big(p^\prime(\rho) \rho^{-1} \div \,(\dot{\rho}\,  W)\big)
-\dot{e} \,\pa_a \big(p^\prime(\rho) \rho^{-1} \div\, (\rho \, W)\big)\tag 3.25 
$$
Since $\dot{\rho}\big|_{\pa\Omega}=0$ and $\dot{e}=\dot{\rho}/\rho$ it follows that
$$\multline
\langle U,\rho\, C^\prime W\rangle\\=\int_{\Omega}\! \Big(\dot{\rho}\, p^{\prime\prime}
\div\,(\rho_{\,}U)\div\,(\rho_{\,}W)
+ p^{\prime}\big(\div\,(\rho_{\,}U)\div\,(\dot{\rho}_{\,}W)
+ \div\,(\dot{\rho}_{\,}U) \div\,(\rho_{\,}W)\big)\Big)\,\rho^{-1}\kappa dy\\
+\int_{\pa\Omega}\!\!
 (-\na_N \dot{p}) U_N W_N \, dS .
\endmultline\tag 3.26
$$
It therefore follows that $\dot{E}\leq C\sqrt{E}(\sqrt{E}+\|F\|)$ and hence
with $E_0=\sqrt{E}$ we have
$$
E_0(t)\leq C\Big( E_0(0)+\int_0^t \|F\|\, d\tau\Big). \tag 3.27
$$

\head  4. Decomposition of the linearized equations into an
operator on the divergence free vector fields and an operator on
the orthogonal complement.\endhead
We will now make an orthogonal decomposition:
${\Cal H}=L^2={\Cal H}_0\oplus{\Cal H}_1$ into divergence free
vector fields ${\Cal H}_0$ and gradients of functions in $H_0^1(\Omega)$; ${\Cal H}_1$.
Let us therefore define the orthogonal projection $P$ onto divergence free
vector fields by
$$
PU^a=U^a-g^{ab}\pa_b p_{U},\qquad \triangle p_U=\div U,
\qquad p_U\big|_{\pa\Omega}=0\tag 4.1
$$
(Here $\triangle q=\kappa^{-1}\pa_a\big( \kappa g^{ab}\pa_b q\big)$. )
$P$ is the orthogonal projection in the inner product, see \cite{L1},
$$
\langle U,W\rangle =\int_\Omega g_{ab} U^a W^b\kappa dy.\tag 4.2
$$
Note also that, with
$$
\|W(t)\|_{r,s}=\sum_{k=0}^s\| \hat{D}_t^k W(t,\cdot)\|_{H^r}\tag 4.3
$$
denoting the Sobolev norms for fixed time with space and time differentiation
of order $r$ and $s$, we have
$$
\|PW\|_{r,s}\leq C\| W\|_{r,s},\qquad
\|(I-P)W\|_{r,s}\leq C\|\div W\|_{r-1,s}\tag 4.4
$$
since its just a matter of solving the Dirichlet problem
and commuting through time derivatives, \cite{L1}.

For a function $f$ that vanishes on the boundary define $A_f
W=-P\big( \nabla \big( W^c \pa_c f\big)\big)$; i.e.
$$
{A}_f W^a=-g^{ab}\pa_b\big( (\pa_c f)
W^c-q\big),\qquad\text{where}\quad \triangle \big( (\pa_c f)
W^c-q\big)=0, \qquad q\big|_{\pa\Omega}=0\tag 4.5
$$
If $U$ and $W$ are divergence free then
$$
\langle U,A_f W\rangle =\int_{\pa\Omega} n_a \, U^a (-\pa_c f) W^c\, dS
\tag 4.6
$$
where $n$ is the unit conormal. If $f\big|_{\pa\Omega}=0$ then
$-\pa_c f\big|_{\pa\Omega} =(-\na_N f) n_c$. It follows that $A_f$
is a symmetric operator on divergence free vector fields, and in
particular the normal operator
$$
A=A_h,\tag 4.7
$$
where $h$ is the enthalpy, is positive, i.e. $\langle W,
A W\rangle \geq 0$, since we assumed the physical condition that
$-\na_N h\geq c_0>0$ on the boundary. The normal operator is order
one, by (4.4)
$$
\|AW\|_{r,s}\leq C\|W\|_{r+1,s}.\tag 4.8
$$

The normal operator has certain delicate
commutator properties with vector fields and positivity properties
which were essential for the existence proof in \cite{L1}.
The main difficulty being that it is not elliptic acting
on vector fields with non vanishing curl.
In order to prove existence one had to replace it by
a sequence of bounded operator which uniformly had the same commutator
and positivity properties.

\comment
In the proof of existence for the equation in the divergence free class
we need to replaced the normal operator
the following smoothed out bounded normal operator.
Let $\rho=\rho(d)$ be a smooth out version of the distance
function to the boundary $d(y)=\dist(y,\pa\Omega)$:
$\rho^\prime\geq 0$, $\rho(d)=d$, when $d\leq 1/4$ and
$\rho(d)=1/2$ when $d\geq 3/4$. Then we can alternatively express $A_f$ as
$$
\underline{A}_f W_a=-\pa_a\big( (f/\rho)(\pa_c \rho) W^c-q\big),\qquad
\triangle \big( (f/\rho)(\pa_c \rho) W^c-q\big)=0,
\qquad q\big|_{\pa\Omega}=0\tag 4.9
$$
Let $\chi(\rho)$ be a smooth function such that $\chi^\prime\geq 0$,
$\chi(\rho)=0$ when $\rho\leq 1/4$, $\chi(\rho)=1$ when $\rho\geq 3/4$.
$A_f$ is unbounded so we now define an approximation that
is a bounded operator: $A_f^\varepsilon W^a=g^{ab} A_f^\varepsilon W_b$,
where
$$
\underline{A}_f^\varepsilon W_a=-\chi_\varepsilon\pa_a\big(
(f/\rho)(\pa_c \rho) W^c\big) +\pa_a q ,\quad
\triangle q
=\kappa^{-1}\pa_a \big( g^{ab} \kappa\chi_\varepsilon
\pa_b\big( (f/\rho)(\pa_c \rho) W^c\big)\big),
\quad q\big|_{\pa\Omega}\!=0\tag 4.10
$$
where $\chi_\varepsilon(\rho)=\chi(\rho/\varepsilon)$. Let $A_f^0=A_f$.
We have
$$
\langle U,A_f^\varepsilon W\rangle
=\int_{\Omega} (f/\rho) \chi_\varepsilon^\prime (\pa_a\rho) \, U^a (\pa_c \rho) W^c\,
\kappa dy
\tag 4.11
$$
from which it follows that $A_f^\varepsilon$ is also symmetric.
And in particular $A^\varepsilon=A^\varepsilon_p$ is positive since we
assumed that $p\geq 0$.
We have
$$
|\langle U,A^\varepsilon_f W\rangle|\leq
\| f/p\|_{L^\infty(\Omega\setminus\Omega_\varepsilon)}\langle U, A^\varepsilon U\rangle^{1/2}
\langle W, A^\varepsilon W\rangle^{1/2}\tag 4.12
$$
where $\Omega_\varepsilon=\{y\in\Omega; d(y,\pa\Omega)<\varepsilon\}$.
Furthermore
$$
\|A_f^\varepsilon W\|_{r,s}\leq C_\varepsilon \|W\|_{r,s}\tag 4.13
$$
\endcomment

We now make the decomposition
$$
W=W_0+W_1,\qquad W_0=P W\in{\Cal H}_0,\qquad W_1=(I-P)W\in
{\Cal H}_1\tag 4.9
$$
\comment 
i.e.
$$
W_1=\nabla q_1,\qquad\triangle q_1=\div W,\quad
q_1\big|_{\pa\Omega}=0.\tag 4.10
$$
\endcomment 
We want to decompose the linearized operator
$$
L W=\ddot{W}+{C} W-B(W,\dot{W}),\qquad
\quad\div W\big|_{\pa\Omega}=0\tag 4.10
$$
where $B$ is a bounded operator and 
$$
{C} W^a=-g^{ab}\pa_b\big( (\pa_c h) W^c+p^{\,\prime}\div
W\big),\tag 4.11
$$
into an operator onto the divergence free part and an operator on the complement.
\comment
Taking the divergence on (4.10) gives a wave equation with Dirichlet boundary
condition for the divergence:
$$
\hat{D}_t^2 \div W-\triangle\big( p^{\,\prime} \div W\big)=
\triangle \big( W^c \pa_c h\big) +\div F+\div B, \qquad\quad\div
W\big|_{\pa\Omega}=0. \tag 4.13
$$
\endcomment 

The projection to highest order commutes with time
differentiation:
$$
P\ddot{W}_0=\ddot{W}_0,\qquad P\ddot{W}_1=PB_2(W_1,\dot{W}_1),\qquad
B_2(W,\dot{W})^a=-g^{ab}\big(\ddot{g}_{bc} W^c+2\dot{g}_{bc}\dot{W}^c\big)\tag 4.12
$$
where $\dot{g}_{ab}\!=\!\check{D}_t g_{ab}$,
 $\ddot{g}_{ab}\!=\!\check{D}_t^2 g_{ab}$, and $\check{D}_t\!=\!D_t\!-\!(\div V)$.
In fact, applying $D_t^2$ to
 $g_{ab} W_1^b\!=\!\pa_q q_1$ gives  $ g_{ab}\ddot{W}_1^b
+2 \dot{g}_{ab}\dot{W}^b+\ddot{g}_{ab} W^b=\pa_a \ddot{q}_1 $.
Here $\ddot{q}_1=D_t^2 q_1$ vanishes on the boundary since $q_1$ does.
The projection of $g^{ab}\pa_b \ddot{q}_1$ therefore vanishes and
(4.12) follows since $\hat{D}_t$ preserves the divergence free
condition. Furthermore with $A$ given by (4.5)-(4.7) and $C$ by (4.11) we have  
$$
P{C} W= A W, \qquad \text{if}\qquad \div W
\big|_{\pa\Omega}=0,\tag 4.13
$$
since the projection of the highest order term,
$\nabla\big(p^{\,\prime}\div W\big)$, vanishes since $\div W
\big|_{\pa\Omega}=0$.

We now want to project (4.10) onto the divergence free vector fields 
using (4.12)-(4.13).
We get
$$
PL W= \ddot{W}_0+A W_0
 +PB_2(W_1,\dot{W}_1)+PAW_1-PB(W,\dot{W})\tag 4.14
$$
where $A$ is the normal operator (4.7). Similarly, applying $(I-P)$ to (4.10) gives
$$
(I-P) LW=\ddot{W}_1-PB_2(W_1,\dot{W}_1)-\nabla \big(p^{\,\prime}\div
W_1\big)- (I-P)\nabla \big( W^c\pa_c h\big) -(I-P)B(W,\dot{W}),\tag
4.15
$$
subject to the boundary condition $\div W_1\big|_{\pa\Omega}=0$.
\comment 
Taking the divergence of (4.16) gives the equation (4.13) for $\div W_1$,
and $W_1$ can alternatively be defined as the solution of (4.10) and (4.13).
\endcomment

\proclaim{Lemma {4.}1}
Let $\tilde{L}$ be defined by
$$\align 
\tilde{L} W_0&= \ddot{W}_0+A W_0,\tag  4.16\\
\tilde{L} W_1&=\ddot{W}_1-PB_2(W_1,\dot{W}_1)-\nabla
\big(p^{\,\prime}\div
W_1\big)\tag 4.17
\endalign 
$$
and let $\tilde{M}$ be defined by 
$$\align 
P\tilde{M}W&=PB_2(W_1,\dot{W}_1)+PAW_1-PB(W,\dot{W}), \tag 4.18\\
(I-P)\tilde{M}W&=-(I-P)\nabla \big( W^c\pa_c h\big)-(I-P)B(W,\dot{W}) \tag 4.19
\endalign 
$$
We have 
$$
LW=\tilde{L} W+\tilde{M}W\tag 4.20
$$
\endproclaim
If $P_0=P$,
$P_1=(I-P)$  and
 ${L}_{ij}=P_i{L}P_j$ then
 $\tilde{L}$ respectively $\tilde{M}$ are essentially the diagonal respectively
the off-diagonal part of ${L}$. 
It turns out that we can invert (4.16) on ${\Cal H}_0$, see section 5, and 
(4.17) on ${\Cal H}_1$, see section 6.
 The interaction term $\tilde{M}$ is lower order,
but in subtle way, since it contains
space derivatives $\pa W$. The estimates for (4.17) gives us 
control of an additional space
derivative of $W_1$ and that is all that is needed to
estimate (4.18). (4.19) also contains a space derivative of
$W_0$. However, in our estimates for (4.17) we can replace
this space derivative by a time derivative and the estimates for (4.16)
gives us control of an additional time derivative. 
The estimates for (4.16)-(4.17) will be summarized in section 7.

\head 5. Existence and estimates in the divergence free class.\endhead
In \cite{L1} we proved existence of solutions for
$$
\ddot{W}_0+AW_0=F_0, \qquad W_0\big|_{t=0}=\tilde{W}_{00},\qquad
\dot{W}_0\big|_{t=0}=\tilde{W}_{10}\tag 5.1
$$
We have
\proclaim{Proposition {5.}1}  Suppose that
$x, h$ are smooth, $h\, \big|_{\pa\Omega}=0$ and
$\na_N h\, \big|_{\pa\Omega}\leq -c_0<0$
for $0\leq t\leq T$.
Then if initial data and the inhomogeneous term in (5.1)
are smooth and divergence free it follows that (5.1) has
a smooth solution for $0\leq t\leq T$.
Furthermore, with a constant $C$ depending only on the norm of $x$ and $h$,
$T$ and the constant $c_0$ we have
$$
E_{0r}(t)\leq
C\Big( E_{0r}(0)+\int_0^t\|F_0(\tau)\|_{H^r}\, d\tau\Big),
\qquad E_{0r}(t)=\|\dot{W}_0(t)\|_{H^r} +\|W_0(t)\|_{H^r}
+\langle W_0(t)\rangle_{r}\tag 5.2
$$
where
$$
\|W(t,r)\|_{H^r}=\|W(t,\cdot\,)\|_{H^r(\Omega)},\qquad \langle
W(t)\rangle_{r}=\|W_N(t,\cdot\,)\|_{H^r(\pa\Omega)}\tag 5.3
$$
and $W_N=N_a W^a$ is the normal component.
\endproclaim
\demo{Proof}
In case $\div V=0$ this was proven in \cite{L1} and
the proof there can be easily modified by multiplying or dividing
by $\kappa\!=\!\det{(\pa x /\pa y)}$.
Let us now indicate what needs to be changed
in \cite{L1} in order to deal with the case $\div V\!\!\neq\! 0$.
We can use the same set of tangential
vector fields as in \cite{L1}, but they are no longer divergence
free so the Lie derivative with respect to these no longer 
preserves the divergence free condition. But one can easily modify the Lie
derivative so it preserves the divergence free condition.
The modified Lie derivative with respect to a vector field $T$ applied to a vector field
$W$ is 
$$
\hat{\Cal L}_T W={\Cal L}_T W+(\div T) W,\tag 5.4
$$
where ${\Cal L}_T W$ is the Lie derivative. 
It satisfies $\div\hat{\Cal L}_T W=\hat{T}\div W$, where for a function 
$f$, $\hat{T}f =Tf+(\div T)f$.
 One then has to make it so
one always apply this modified Lie derivative to vector fields.
However we use the usual Lie derivative, when applied to one forms
since it commutes with covariant differentiation.
In deriving the estimates for all components of a vector field in terms
of the divergence, the curl and the tangential components we
use ${\Cal L}_T (g_{ab} W^b)= g_{ab}\hat{\Cal L}_T W^b
+(\check{\Cal L}_T g_{ab}) W^b$, where  $\check{\Cal L}_T g_{ab}
={\Cal L}_T g_{ab}-(\div T) g_{ab}$. Let us now examine how the
critical commutator with the normal operator, (4.5), is changed from
what it was in \cite{L1}. With $A_f$ given by (4.5) and $\underline{A}_f W_a=
g_{ab} A_f W^b$ we have 
$$
{\Cal L}_T \underline{A}_f W_a={\Cal L}_T \pa_a \big((\pa_c f) W^c-q\big)
=\pa_a \big( (\pa_c f) \hat{\Cal L}_T W+\pa_c (\check{T} f) W^c- Tq
+f(\pa_c \div T) W^c\big)\tag 5.5
$$
 where $\check{T}f=Tf-(\div T) f$. When we project again the last two terms
 vanish since they vanish on the boundary so the commutator relation in
\cite{L1} will be replaced by
$$
P{\Cal L}_T \underline{A}_f W =\underline{A}_f \hat{\Cal L}_T W+
\underline{A}_{\check{T}f} W\tag 5.6
$$
The issue of how to deal with the initial conditions in case $\div
V\neq 0$ was discussed in section 3. \qed\enddemo

Now, the norms used in Proposition {5.}1 are natural for the initial
value problem. However, when solving the wave equation with
Dirichlet boundary conditions it is more natural to first look
on norms with many time derivatives. Because of the coupling
between the two equations we must therefore also estimate more time
derivatives of the divergence free part.
From differentiating (5.1) we get:

\proclaim{Lemma {5.}2} Suppose that $W$ is a smooth solution of (5.1) for
$0\leq t\leq T$. Let $E_{0,r}$ be as in (5.2) and
$$
{\|}W(t){\|}_r=\sum_{s+k\leq r}\|W(t)\|_{k,s}
\qquad \text{where}\qquad
\|W(t)\|_{r,s}=\sum_{k=0}^s \| \hat{D}_t^k W(t)\|_{H^r}\tag 5.7
$$
Then for $0\leq t\leq T$ we have
$$
{\|}\dot{W}_0{\|}_r+{\|}W_0{\|}_r\leq C\big(E_{0,r}+{\|} F_0{\|}_{r-1}\big),
\qquad r\geq 1
\tag 5.8
$$
\endproclaim
\demo{Proof} The proof is just differentiation of (5.1) using that
$A$ is order one, (4.8);
$$
{\|}\ddot{W}_0{\|}_{r-1}\leq C({\|}W_0{\|}_r+{\|}F_0{\|}_{r-1}),
\qquad r\geq 1\tag 5.9
$$
which proves (5.8) for $r=1$ so we may assume that $r\geq 2$ in (5.8).
We have
$$
{\|} W_0{\|}_r\leq \|W_0\|_{r,0}+\|\dot{W}_0\|_{r-1,0}+{\|}\ddot{W}_0{\|}_{r-2},
\qquad r\geq 2
\tag 5.10
$$
so together with (5.9) we get
$$
{\|} {W}_0{\|}_{r}\leq C\big(E_{0,r}+{\|} {W}_0{\|}_{r-1}
+{\|} F_0{\|}_{r-2}\big),\qquad r\geq 2\tag 5.11
$$
Since also ${\|}W_0{\|}_1\leq E_{0,1}$ we can use induction
in $r$ to prove that
$$
{\|}W_0{\|}_r\leq C\big(E_{0,r}+{\|} F_0{\|}_{r-2}\big),\qquad
r\geq 2\tag 5.12
$$
Similarly
$$
{\|}\dot{W}_0{\|}_r\leq \|\dot{W}_0\|_{r,0}+{\|}\ddot{W}_0{\|}_{r-1},
\qquad r\geq 1\tag 5.13
$$
so by (5.9) again
$$
{\|}\dot{W}_0{\|}_r\leq C\big(E_{0,r}+{\|}W_0{\|}_r+{\|} F_0{\|}_{r-1}\big),
\qquad r\geq 1\tag 5.14
$$
(5.8) for $r\geq 2$ now follows from (5.12) and (5.14).
\qed\enddemo

\proclaim{Theorem {5.}3} With notation and assumptions as in
Proposition {5.}1 and Lemma {5.}2 we have
$$
\tilde{E}_{0r}(t)\leq C\Big( \tilde{E}_{0r}(0)+\int_0^t{\|}
F_0(\tau){\|}_{r}\, d\tau\Big), \qquad
\tilde{E}_{0r}(t)={\|}\dot{W}_0(t){\|}_{r} +{\|}W_0(t){\|}_{r}
+\langle W_0(t)\rangle_{r}\tag 5.15
$$
\endproclaim
\demo{Proof} (5.15) follows from Proposition {5.}1 since
$$
{\|} F_0(t){\|}_{r-1}\leq {\|} F_0(0){\|}_{r-1}+\int_0^t {\|}
\dot{F}_0(\tau){\|}_{r-1}\, d\tau , \qquad {\|} F_0(0){\|}_{r-1}\leq
C\big({\|} \dot{W}(0){\|}_{r}+{\|} W(0){\|}_{r}\big)\qed\tag 5.16
$$
\enddemo

\head 6. Existence and estimates for the wave equation.\endhead

We consider the Cauchy problem for the wave equation on a bounded domain
with Dirichlet boundary conditions:
$$\align
\hat{D}_t^2(e^\prime \psi)-\triangle \psi &=f,\qquad
\text{in}\quad [0,T]\times\Omega,\qquad \psi\big|_{\pa\Omega}=0,\tag 6.1 \\
\psi\big|_{t=0}\!\! &=\tilde{\psi}_0,\quad
{D}_t\psi\big|_{t=0}\!\!=\tilde{\psi}_1 \tag 6.2
\endalign
$$
Here
$$
\triangle\psi =\frac{1}{\sqrt{\det{g}}}\pa_a
\Big( \sqrt{\det{g}} g^{ab}\pa_b \psi\Big),\tag 6.3
$$
where $g^{ab}$ is the inverse of the metric $g_{ab}$ and
$\det{g}=\det\{g_{ab}\}=\kappa^2$, in our earlier notation.
We assume that  $g^{ab}$ is symmetric (since the metric is), and that $g^{ab}$
and $e^\prime$ are smooth satisfying:
$$
0<e^\prime+1/e^\prime<c_1^\prime,\qquad 
\sum_{a,b=1}^{n} \big(| g^{ab}|+|g_{ab}|\big)\leq n \,c_1^2
\tag 6.4
$$
for some constants $0<c_1<c_1^\prime<\infty$.

Existence of solutions for a wave equation with Dirichlet
boundary conditions and initial conditions satisfying
some compatibility conditions is well known, see e.g. \cite{H,Ev}.
In order for (6.1)-(6.2) to be solvable initial data must be compatible
with the boundary condition. If we move the Laplacian in (6.1)
over to the right hand side and differentiate (6.1)
with respect to time we get
$$
{D}_t^{k+2} \psi=b_k\big(\psi,...,{D}_t^{k+1} \psi,
\pa \psi,...,\pa{D}_t^{k} \psi,\pa^2\psi,...,
\pa^2 {D}_t^{k} \psi\big)+D_t^k f\tag 6.5
$$
for some functions $b_k$. We therefore define functions of the
space variables only $\tilde{\psi}_k$
$$
\tilde{\psi}_{k+2}=b_k\big(\tilde{\psi}_0,...,\tilde{\psi}_{k+1},\pa\tilde{\psi}_0,...,\pa
\tilde{\psi}_{k},
\pa^2\tilde{\psi}_0,...,\pa^2\tilde{\psi}_{k}\big)\big|_{t=0}+D_t^k
f\big|_{t=0}\tag 6.6
$$
where $\tilde{\psi}_0$ and $\tilde{\psi}_1$ are as in (6.2). For
this to be compatible with the boundary conditions we must have
$$
\tilde{\psi}_k\big|_{\pa\Omega}=0,\qquad \text{for}\quad k\leq
m-1\tag 6.7
$$
(6.7) is called the $m^{th}$ order compatibility condition. Since
$\tilde{\psi}_k$ are determined from the initial conditions
$\tilde{\psi}_0$ and $\tilde{\psi}_1$ this gives some
compatibility conditions on the initial conditions. We have:

\proclaim{Proposition {6.}0} Suppose that $g,e^\prime$ are smooth
satisfying (6.4). Then if initial data
$(\tilde{\psi}_0,\tilde{\psi}_1)$ and $f$ are smooth and satisfy
the $m^{th}$ order compatibility condition for all $m$, it follows
that (6.1)-(6.2) has a smooth solution $\psi$.
\endproclaim
\demo{Proof}
The result in \cite{H} is stated with
vanishing initial conditions but, if the compatibility conditions are
satisfied, one can reduce to that case by subtracting off an approximate
solution satisfying the equation to all orders as $t\to 0$.
Let
$$
\tilde{\psi}=\sum_{k=0}^{\infty}\chi(t/\varepsilon_k)
t^k\tilde{\psi}_k /k!,\tag 6.8
$$
where $\chi$ is smooth $\chi(s)=1$ for $|s|\leq 1/2$ and $\chi(s)=0$
for $|s|\geq 1$, and the
 sequence $\varepsilon_k>0$ are chosen small enough so that the series
converges in $C^m([0,T],H^m)$ for any $m$. This is obtained if
take $(\|\tilde{\psi}_k\|_{H^k}+1)\varepsilon_k\leq 1/2$. Then
$\overline{\psi}=\psi-\tilde{\psi}$ satisfies (6.1) with vanishing
initial conditions and a right hand side $\overline{f}$ that
vanishes to all orders as $t\to 0$:
$$
\hat{D}_t^2(e^\prime \overline{\psi})-\triangle\overline{\psi}=
-\hat{D}_t^2(e^\prime\tilde{\psi})+\triangle\tilde{\psi}+f=\overline{f}
\tag 6.9
$$
For this case existence of a smooth solution $\overline{\psi}$ to
(6.9) follows from Theorem 24.1.1 in \cite{H}. Since the theorem
in \cite{H} is more general let
us just point out the main steps needed for our case. Existence
follows from duality, using the Hahn-Banach extension theorem and
the Riesz representation theorem. For this one has to show estimates
for the adjoint operator in negative Sobolev spaces. Suppose that
$\psi$ satisfy (6.1) and let ${\psi}_N=(I-\triangle)^{-N}\psi$,
where $N\geq 0$ and $\triangle$ is the Dirichlet Laplacian, i.e.
inductively, we define ${\psi}_k$ to be the solutions of
$(I-\triangle)\psi_{k+1}=\psi_k$, with boundary conditions
$\psi_{k+1}\big|_{\pa\Omega}=0$. Then $\psi_N$ satisfies (6.1)
with $f$ replaced by $f_N
+(I-\triangle)^{-N}[\hat{D}_t^2,(I-\triangle)^{N}]\psi_N$, where
$f_N= (I-\triangle)^{-N} f$. The norm of this is bounded by
$\|f_N\|+\|\psi_N\|+\|D_{t\,}\psi_N\|$.  Using the energy estimate
in Lemma {6.}2 then gives us an estimate $\|D_t
\psi_N(t,\cdot)\|+\|\nabla\psi_N(t,\cdot)\| \leq C\int_0^t
\|f_N\|\, d\tau $ \qed\enddemo

\proclaim{Lemma {6.}2}
 Suppose that $g^{ab},e^\prime,f $ and $\psi$ are smooth and satisfy
(6.1)-(6.4) for $0\leq t\leq T$. 
Let $\nabla^a=g^{ab}\pa_b$ and, for $r\geq 1$,
$$
E(t)=\Big(\sum_{s=0}^{r-1} \frac{1}{2}\!\int_\Omega{ \!\! \big( e^\prime
(\!{D}_t^{\,s+1}\psi)^2\! + |\hat{D}_t^{s}\nabla \psi|^2+\psi^2\, \big)
\kappa d y}\Big)^{1/2}\!\!\!\!\!\!, \tag 6.10
$$
Then
$$
\frac{d E}{dt} \leq C \big( E+\|f\|_{0,r-1}\big),
\qquad\text{where}\quad
 \| \phi\|_{r,s}=\sum_{k\leq s,\, |\alpha|\leq r} \|D_t^k \pa_y^\alpha \phi\|
\tag 6.11
$$
\endproclaim
\demo{Proof} We will prove that $ d E^2/dt$ is bounded by $E$
times the right hand side of (6.11) and (6.11) follows from this
since $d E/dt=(d E^2/dt)/(2E)$.
Since  $D_t\kappa=\kappa \div V$, we have with
 $\hat{D}_t=D_t+\div V$ and $\check{D}_t=D_t-\div V$:
$$
\multline
\frac{d E^2}{dt}=\sum_{s\leq r-1}
\int_\Omega{ \Big( e^\prime ({D}_t^{s+1}\psi)({D}_t^{s+2}\psi)+
g_{ab}\big(\hat{D}_t^s\nabla^a \psi\big)
\big(\hat{D}_t^{s+1}\nabla^b \psi\big)
\Big)\,  \kappa \,d y}\\
+\frac{1}{2}\int_\Omega{ \Big((\hat{D}_t e^\prime) ({D}_t^{s+1}\psi)^2+
(\check{D}_t g_{ab})\big(\hat{D}_t^s\nabla^a \psi\big)
\big(\hat{D}_t^{s}\nabla^b \psi\big)
+\kappa^{-1} (D_t\kappa)\psi^2 +2\psi D_t\psi \Big)\,  \kappa \,d y}.
\endmultline\tag 6.12
$$
Here the terms on the second row are bounded by a constant times $E^2$.
Applying $D_t^{s+1}$ to $\pa_a \psi=g_{ab}\nabla^b\psi$ gives 
$\pa_a D_t^{s+1} \psi=\sum_{i=0}^{s+1} \binom{s+1}{i} 
(\check{D}_t^{s+1-i} g_{ab})\hat{D}_t^{i} \nabla^b \phi$ so 
 $$
 \hat{D}_t^{s+1}\nabla^a\psi =g^{ab}\pa_a
 D_t^{s+1}\psi-
 \sum_{i=0}^{s}{\tsize{\binom{s+1}{i}}} g^{ab} (\check{D}_t^{s+1-i} g_{bc}) \hat{D}_t^i
 \nabla^c\psi.\tag 6.13
 $$
Up to terms bounded by a constant times $E^2$,
 (6.12) is therefore equal to 
$$
\multline
\sum_{s\leq r-1}
\int_\Omega{ \!\Big( e^\prime ({D}_t^{s+1}\psi)({D}_t^{s+2}\psi)+
\big(\hat{D}_t^s\nabla^a \psi\big)
\big(\pa_a {D}_t^{s+1}\psi\big)
\Big)\,  \kappa \,d y}\\
=\sum_{s\leq r}\int_\Omega{ \Big( ({D}_t^{s+1}\psi )
\big(e^\prime{D}_t^{s+2}\psi-\kappa^{-1}
\pa_a\big(\kappa \hat{D}_t^s \nabla^a \psi\big)
\Big)\, \kappa dy}
\endmultline\tag 6.14
$$
where we have integrated  by parts. If we apply ${\hat D}_t^{s}$ to 
$\hat{D}_t^2\big(e^\prime \psi\big)-\kappa^{-1}\pa_a\big(\kappa \nabla^a \psi\big)
=f$ we obtain 
$$
\hat{D}_t^s \Big( \hat{D}_t^2\big( e^\prime \psi\big)
 -\kappa^{-1}\pa_a\big(\kappa \nabla^a
\psi\big)\Big)=
e^\prime {D}_t^{s+2}\psi -\kappa^{-1}\pa_a\big(\kappa \hat{D}_t^s\nabla^a
\psi\big)+\sum_{i=0}^{s+1} {\tsize{\binom{s+2}{i}}} (\hat{D}_t^{s+2-i}
 e^\prime)({D}_t^{i}\psi). 
\tag 6.15
$$
Since the $L^2$ norm of the last term is bounded
by $CE$ plus the $L^2$ norm of $\psi$ the lemma follows.
\qed\enddemo

One can get additional space
regularity from taking time derivatives of the equation (6.1) and solving
the Dirichlet problem for the Laplacian.

\proclaim{Lemma {6.3}}  Suppose that $g^{ab}$ and $e^\prime$ are smooth and satisfy
(6.4) and that $f$ is smooth, for $0\leq t\leq T$.
Suppose also that $\psi$ is a smooth solution of (6.1), for $0\leq t\leq T$.
Let $\|\psi\|_{s,r}$ be as in Theorem {6.}1 and let
$
{\|} \psi{\|}_{r}=\sum_{s+k\leq r} \|D_t^s \psi\|_k.
$
Then
$$
{\|} \psi{\|}_r\leq C\big(\|\psi\|_{0,r}+\|\psi\|_{1,r-1}+{\|}f{\|}_{r-2}\big).\tag 6.16
$$
\endproclaim
\demo{Proof} Since 
$\triangle \psi=\hat{D}_t^2\big(e^\prime\psi\big)-f$ and 
$ \psi\big|_{\pa\Omega}=0$
it follows that 
$$
\triangle {D}_t^s \psi=\sum_{i=0}^{s+2} {\tsize{\binom{s+2}{i}}}
 (\hat{D}_t^{s+2-i}e^\prime) {D}_t^{i} \psi
-\hat{D}_t^s f-[\hat{D}_t^s\triangle -\triangle D_t^s]\psi
\qquad {D}_t^s\psi\big|_{\pa\Omega}=0\tag 6.17
$$
so by the standard elliptic estimates 
$$
\|{D}_t^s \psi\|_{H^{k+2}}
\leq C\big( \sum_{i=0}^{s+2}\|{D}_t^{i}\psi\|_{H^{k}}
+\|\hat{D}_t^s f\|_{H^{k}}
+\sum_{i=0}^{s-1}\|{D}_t^i\psi\|_{H^{k+2}}\big)\tag 6.18
$$
Here the last term is lower order and is absent if $s=0$ so using
induction in $s$ we get
$$
\sum_{i=0}^s\|{D}_t^i \psi\|_{H^{k+2}}
\leq C \sum_{i=0}^{s}\big(\|{D}_t^{i+2}\psi\|_{H^k}
+\|{D}_t^i f\|_{H^k}\big)\tag 6.19
$$
or with $\|\psi\|_{r,s}=\sum_{k=0}^s\|D_t^k \psi\|_{H^r}$
$$
\|\psi\|_{s+2,r-s-2}\leq C\big(\|\psi\|_{s,r-s}+\|f\|_{s,r-s-2}\big),
\qquad 0\leq s\leq r-2\tag 6.20
$$
Since ${\|}\psi{\|}_r=\sum_{s=0}^r \|\psi\|_{r-s,s}$ it therefore
inductively follows that
$$
{\|}\psi{\|}_r\leq C\big(\|\psi\|_{0,r}+\|\psi\|_{1,r-1}+{\|}f{\|}_{r-2}
\big)\qed \tag 6.21
$$
\enddemo

Summing up;
\proclaim{Proposition {6.}4} There are constants $C_r$ such that the solution of 
(6.1) satisfy 
$$
{\|}\psi(t,\cdot){\|}_r\leq C_r\Big({\|}\psi(0,\cdot){\|}_r
+\int_0^t{\|}\dot{f}{\|}_{r-2}\, d\tau\Big),\qquad r\geq 2\tag 6.22 
$$
\endproclaim
\demo{Proof} It follows from Lemma {6.}2 that 
$E(t)\leq C_r \big( E(0)+\int_0^t\|f\|_{0,r-1}\, d\tau \big) $.
 Using (6.13) we see that the energy energy $E$ in (6.10)
is equivalent to $\|\psi\|_{0,r}+\|\psi\|_{1,r-1}$. Furthermore 
${\|}f(t,\cdot){\|}_{r-2}\leq {\|}f(0,\cdot){\|}_{r-2}
+\int_0^t {\|}\dot{f}(t,\cdot){\|}_{r-2}\, d\tau $ and
since also ${\|}f{\|}_{r-2}\leq C{\|}\psi{\|}_r$ the proposition follows from 
Lemma {6.}3. 
\qed\enddemo

As pointed out in section 4 we actually want to solve the
equation:
$$\align
\ddot{W}_1-PB_2(W_1,\dot{W}_1)-\nabla \big(p^{\,\prime}\div
W_1\big)&= F_1,\qquad \div W_1\big|_{\pa \Omega}=0 \tag
6.23\\
W_1=\nabla q_1,\qquad q_1\big|_{\pa\Omega}=0,\tag 6.24
\endalign
$$
where $(I-P)F_1=F_1$, 
which is equivalent to 
$$\align 
\hat{D}_t^2 \phi-\triangle\big(p^{\,\prime} \phi
\big)&=\div F_1,\qquad \phi|_{\pa\Omega}=0\tag 6.25\\
W_1=\nabla q_1,\qquad 
\triangle q_1&=\phi,\qquad q_1\big|_{\pa\Omega}=0. \tag 6.26 
\endalign 
$$
In fact, for $W_1$ of the form (6.24) the left hand side of (6.23) is
$(I-P)\ddot{W}_1-(I-P)\nabla \big(p^\prime \div W_1\big)$, and 
$(I-P)H=0$ is equivalent to $\div H=0$. Assuming that the compatibility conditions 
are satisfied we can solve (6.25)-(6.26) and this then also gives us a solution of
(6.23)-(6.24). The initial conditions for (6.25) are
$$
\phi\big|_{t=0}=\tilde{\phi}_0,\qquad \hat{D}_t \phi\big|_{t=0}=\tilde{\phi}_1
\tag 6.27
$$

\proclaim{Lemma {6.}5} Suppose that $W_1$ satisfies (6.23) and set 
$$
E(t)=\Big(\frac{1}{2}\sum_{s=0}^{r}
\int_{\Omega} \big( |\hat{D}_t^{s+1} {W}_1^{}|^2 
+p^\prime |\div \,( \hat{D}_t^s W_1^{})|^2+|W_1|^2\big)
\kappa \, dy\Big)^{1/2}.\tag 6.28 
$$
Then 
$$
\frac{d E}{dt}\leq C\big( E+\|F_1\|_{0,r}\big).\tag 6.29 
$$
\endproclaim
\demo{Proof} Let $W_{1k}=\hat{D}_t^k W_1$, 
$$\multline 
\frac{d E^2}{dt}=\sum_{s=0}^{r} \int_{\Omega}
\big( \dot{W}_{1s} \cdot \ddot{W}_{1s}+ p^\prime \div W_{1s} \div \dot{W}_{1s}\big)\, 
\kappa dy\\
+\frac{1}{2}\int_{\Omega}\big( (\check{D}_t g_{ab})\dot{W}_{1s}^a\dot{W}_{1s}^b
+(\check{D}_t p^\prime)( \div W_{1s})^2 +\kappa^{-1} (D_t\kappa ) |W_1|^2 
+2\, W_1\cdot \dot{W}_1\big) \kappa dy .
\endmultline  \tag 6.30 
$$
Integrating by parts, using that $\div W_{1s}\big|_{\pa\Omega}=0$,
it therefore follows that 
$$
\frac{d E^2}{dt}\leq \sum_{s=0}^{r} \int_{\Omega}  \dot{W}_{1s}\cdot 
\big(\ddot{W}_{1s}-\nabla\big(p^\prime \div W_{1s}\big)\big)\, \kappa dy+CE^2 . \tag 6.31 
$$
Using (6.23) this proves (6.29) for $r=0$. To prove it for $r\geq 1$ we have to commute
time derivatives through (6.23), which can be written 
$$
g_{ab}\ddot{W}_1^b -\pa_a\big( p^\prime \div W_1\big)=g_{ab}\tilde{B}_0^b, 
\qquad
\tilde{B}_0(W_1,\dot{W}_1,F_1)=F_1+P B_2(W_1,\dot{W}_1)\tag 6.32 
$$
With $q=p^\prime/\kappa$ we have 
$$\multline 
D_t^s \big( q^{-1} \pa_a\big(p^\prime \div W_1\big) \big)
= D_t^s \pa_a \big(\kappa \div W_1\big)
+D_t^s \big( (\pa_a \ln q)
\kappa \div W_1\big)\\
= \pa_a \big(\kappa \div W_{1s}\big)
+\sum_{k=0}^s{\tsize{\binom{s}{k}}} (\pa_a D_t^{s-k} \ln q)\kappa \div W_{1k}
=q^{-1} \pa_a\big(p^\prime\div W_{1s}\big) 
+\sum_{k=0}^{s-1}{\tsize{\binom{s}{k}}} (\pa_a D_t^{s-k} \ln q)\kappa \div W_{1k}
\endmultline \tag 6.33 
$$
Hence 
$$
q D_t^s \big( q^{-1} \pa_a\big(p^\prime \div W_1\big) \big)
= \pa_a\big(p^\prime\div W_{1s}\big) +
\sum_{k=0}^{s-1}{\tsize{\binom{s}{k}}} (\pa_a D_t^{s-k} \ln q)p^\prime \div W_{1k}
\tag 6.34
$$
Multiplying (6.32) by $q^{-1}$, applying $D_t^s$ and dividing by $q^{-1}$ therefore gives 
$$
g_{ab}\ddot{W}_{1s}^b -\pa_a\big( p^\prime \div W_{1s}\big)=
g_{ab} \tilde{B}_s^b\big(W_1,...,\dot{W}_{1s},\div W_1,...,\div W_{1s-1}, 
F_1,...,\hat{D}_t^s F_1\big)\tag 6.35
$$
where $\tilde{B}_s$ is a bounded operator of its arguments:
$$
\|\tilde{B}_k\|\leq C\sum_{k=0}^s \big(\|W_{1k}\|+\|\dot{W}_{1k}\|+
\|\hat{D}_t^k F_1\|+\|\div W_{1k}\|\big) \tag 6.36
$$
This together with (6.31) proves (6.29) also for $r\geq 1$. 
\qed\enddemo

 \proclaim{Theorem {6.}6} Suppose that the initial conditions 
(6.27) and the inhomogeneous term in (6.25) are smooth and
satisfy the compatibility conditions for all orders. 
Then (6.25)-(6.27) has a smooth $\phi$. 
Furthermore  with $W_1$ given by (6.26) we have 
 $$
 E_{1r}(t)\leq C_r E_{1r}(0)+C_r\int_0^t {\|} \dot{F}_1(\tau){\|}_{r-1}\,
 d\tau ,\qquad
 E_{1r}(t)={\|} W_1(t){\|}_{r+1}\tag 6.37
 $$
 for $r\geq 1$, and for $r=0$ the same inequality holds with
 ${\|} \dot{F}_1(\tau){\|}_{r-1}$ replaced by
 ${\|} {F}_1(\tau){\|}_{r}$.

\endproclaim
\demo{Proof} First we assume that $r\geq 2$. By the second part of (4.4) 
we see that with $\phi$ and $W_1$ as in (6.37)
$$
{\|}\phi{\|}_r\leq {\|}\pa W_1{\|}_r+{\|}W_1{\|}_r\leq C{\|}\phi{\|}_r\tag 6.38
$$
where $\pa$ stands for space derivatives only. 
Furthermore by (6.23) 
$$
{\|}W_1{\|}_{r+1} \leq {\|}\pa W_1{\|}_r+{\|}W_1{\|}_r
+{\|} \ddot{W}_1{\|}_{r-1}\leq C\big(  {\|}\pa W_1{\|}_r+{\|}W_1{\|}_r
+{\|}F_1{\|}_{r-1}\big)\tag 6.39
$$
Since also 
$$
{\|}F_1(t){\|}_{r-1}\leq {\|}F_1(0){\|}_{r-1}
+\int_0^t {\|}\dot{F}_1{\|}_{r-1} \, d\tau 
\leq C{\|}W_1(0){\|}_{r+1}+\int_0^t {\|}\dot{F}_1{\|}_{r-1} \, d\tau 
\tag 6.40
$$
(6.37) for $r\geq 2$ follows from Proposition {6.}4. 
For $r=1$, (6.37) follows from Lemma {6.}5 and the fact that by the 
second part of (4.4)
$\|\pa D_t W_1\|\leq C\big(\|D_t \div W_1\|+\|\div W_1\|\big)$
and by (6.23) $\|\pa^2 W_1\|\leq C\|\pa \div W_1\|\leq C^\prime\big(\|\ddot{W}_1\|
+\|D_t \pa W_1\|+ {\|}W{\|}_1\big)+\|F_1\|$. 
\qed\enddemo

\comment
We can similarly turn this around and get estimates for
$E_r$ from more space derivatives instead which can be used
for the initial conditions.
\proclaim{Lemma {6.}4} Suppose that $g^{ab}$ and $e^\prime$ are smooth and satisfy
(6.4) and that $f$ is smooth, for $0\leq t\leq T$.
Suppose that that $\psi$ is a smooth solution of (6.1), for $0\leq t\leq T$.
Let $E_{2,r}(t)=\|\psi(t)\|_{r,0}+\|\psi(t)\|_{r-1,1}$. Then
$$
{\|} \psi{\|}_r\leq C\big( E_{2,r}+{\|}f{\|}_{r-2}\big),
\qquad r\geq 2,\qquad
{\|}\psi{\|}_1\leq C E_{2,1}\tag 6.23
$$
\endproclaim
\demo{Proof} The proof is exactly the same as that of Lemma {6.}3
apart from that the Laplacian and the second derivative with respect
to time changed role. We obtain the result simply by
differentiating
$$
\hat{D}_t^2 \psi=(e^\prime)^{-1} \big( \triangle \psi+f)\tag 6.24
$$
with respect to time.
$$
\|\psi\|_{r-k,k}\leq C\Big(\|\psi\|_{r,0}+\|\psi\|_{r-1,1}
+\sum_{j=2}^{k}\|f\|_{r-j,j-2}\Big) \tag 6.25
$$
 This follows from (6.5) using induction in $k$.
\qed\enddemo
\endcomment

\head 7. The proof of the theorem.\endhead 
We are now in position to prove Theorem {3.}1.
We want to show that 
$$
LW=F,\qquad\div W\big|_{\pa\Omega}=0,
\tag 7.1
$$
with initial conditions 
$$
 W\big|_{t=0}=\tilde{W}_0,\qquad  \dot{W}\big|_{t=0}=\tilde{W}_1,\tag 7.2
$$
has a smooth solution $W$ if $\tilde{W}_0,\tilde{W}_1$ and $F$
are smooth and satisfy the compatibility conditions in section 3 to all orders.
If these compatibility conditions hold then
 we can find an approximate solution $\tilde{W}$ satisfying 
the initial conditions and the equation to all orders as $t\to 0$.
Subtracting off this approximate solution reduces it to finding a smooth solution 
to (7.1) when $\tilde{W}_0=\tilde{W}_1=0$ and $F$ vanishes to all orders as $t\to 0$. 

With $\tilde{L}$ and $\tilde{M}$ are
as in Lemma {4.}1 we have reduced to finding a smooth solution of 
$$
\tilde{L}W=\tilde{F},\qquad\div{W}\big|_{\pa\Omega}=0, \tag 7.3
$$
where
$$
\tilde{F}=F-\tilde{M}W.\tag 7.4
$$
That the operator $\tilde{L}$ is invertible follows from using the 
decomposition $W=W_0+W_1$, in Lemma {4.}1 and applying Theorem {5.}3 to the 
divergence free part $W_0$ and Theorem {6.}6 to $W_1$.
This in fact gives estimates for the solution of (7.3) that can be used to 
show existence for (7.3) with $\tilde{F}$ given by (7.4)  by iteration,
and we hence obtain a solution to (7.1). Let us introduce the norms: 
$$
\||W\||_{r,1}={\|}\dot{W}_0{\|}_r+{\|}W_0{\|}_r 
+\langle W_0\rangle_{r}+{\|}W_1{\|}_{r+1},\qquad W_0=PW,\qquad W_1=(I-P)W
\tag 7.5
$$
and
$$
\|| F\||_{r,2}={\|}F_0{\|}_r +{\|}\dot{F}_1{\|}_{r-1},\qquad F_0=PF,
\qquad F_1=(I-P)F. \tag 7.6
$$
It now follows from Lemma {4.}1, Theorem {5.}3 and Theorem {6.}6: 
\proclaim{Theorem {7.}1} Suppose that $(x,h)$ is a smooth solution to (2.7),
(for $0\!\leq \!t\!\leq \!T\!$), such that $h\big|_{\pa\Omega}\!\!=\!0$ and
 $\na_{N} h|_{\pa\Omega}\!\leq\! -c_0\!<\!0$.   
Suppose also that $\tilde{F}$, $\tilde{W}_0$ and
$\tilde{W}_1$ are smooth and such that there is a smooth function $\tilde{W}$
satisfying the initial conditions (7.2), the boundary condition 
$\div\tilde{W}\big|_{\pa\Omega}\!\!=0$ and $\tilde{L} W=F$ to all 
orders as $t\!\to\! 0$, i.e. $D_t^k\big(L\tilde{W}\!\!-\!\tilde{F})\big|_{t=0}\!\!=\!0$, for
$k\!\geq\! 0$ 
Then (7.2)-(7.3) has a smooth 
solution $W$. 

Furthermore, there are constants $C_r$ such that for any smooth 
solution of (7.3) we have 
$$
\|| W(t)\||_{r,1}\leq C_r\Big(\|| W(0)\||_{r,1}+\int_0^t\||\tilde{F}\||_{r,2} \, d\tau\Big),
\qquad r\geq 1 .
\tag 7.7
$$

Moreover, 
$$
\||\tilde{M} W\||_{r,2}\leq C_r\||W\||_{r,1}.\tag 7.8 
$$
\endproclaim
We remark that the compatibility conditions in the theorem 
are in particular true if $\tilde{W}_0=\tilde{W}_1=0$ and 
$D_t^k \tilde{F}\big|_{t=0}=0$, for $k\geq 0$, since then we can take $\tilde{W}=0$.
Therefore if  $D_t^k F\big|_{t=0}=0$, for $k\geq 0$ and we set $W^0=0$,
it follows that we can inductively solve, for $k\geq 0$: 
$$
\tilde{L}W^{k+1}=F-\tilde{M}W^k,\qquad\div{W}^{k+1}\big|_{\pa\Omega}=0,\qquad
W^{k+1}\big|_{t=0}=\dot{W}^{k+1}\big|_{t=0}=0,\tag 7.9
$$
We claim that $W^k$ converges to a solution of (7.3)-(7.4) and hence to (7.1),
in case $D_t^k F\big|_{t=0}=0$, for $k\geq 0$ and $\tilde{W}_0=\tilde{W}_1=0$.
Since we have already reduced solving (7.1)-(7.2) to this case this would prove the
existence part of Theorem {3.}1. That $W^k$ converges to a smooth solution of (7.3)-(7.4)
follows from using the estimate in Theorem {7.}1. In fact,
$\tilde{L}(W^1-W^0)=F$, and for $k\geq 1$
$\tilde{L}({W}^{k+1}-W^k)=-\tilde{M}({W}^k-W^{k-1})$.
It therefore follows from Theorem {7.}1 that 
$$
 M_{N}\leq C_r\int_0^t \Big(\||F\||_{r,2}
+M_{N}\Big)\, d\tau, \qquad\text{where}\qquad 
 M_{N}=\sum_{k=0}^{N} \||{W}^{k+1}-W^k\||_{r,1}\tag 7.10
$$
It now follows from a standard Gr\"onwall type of argument that 
$$
M_N(t)\leq C_r e^{C_r t}\int_0^{t}\||F\||_{r,2}\, d\tau\tag 7.11
$$
for any $N$. It follows from this that $W^k$ converges to a smooth solution 
$W$ of (7.3)-(7.4) for $0\leq t\leq T$
and therefore we have proven existence of smooth solutions for (7.1)-(7.2). 

Haven proven existence of a smooth solution to (7.1)-(7.2) we now 
also need to prove the estimate in Theorem {3.}1. Applying the estimate in Theorem
{7.}1  gives
$$
\|| W(t)\||_{r,1}\leq C_r\Big(\|| W(0)\||_{r,1}+\int_0^t\||{F}\||_{r,2}+
\|| W\||_{r,1} \, d\tau\Big).
\tag 7.12
$$
By the same Gr\"onwall type of argument as above we get 
$$
\|| W(t)\||_{r,1} \leq C_r e^{C_r t}
\Big(\|| W(0)\||_{r,1}+\int_0^{t}\||F\||_{r,2}\, d\tau\Big),\qquad r\geq 1.\tag 7.13
$$
It therefore only remains to observe that the norms in Theorem {3.}1
are equivalent to those here. It follows from the continuity of the projection (4.4) that 
$$
{\|}\dot{W}{\|}_r+{\|}W{\|}_r\leq {\|}\dot{W}_0{\|}_r+{\|}\dot{W}_1{\|}_r
+ {\|}W_0{\|}_r+{\|}W_1{\|}_r\leq C_r\big({\|}\dot{W}{\|}_r+{\|}W{\|}_r\big)\tag 7.14
$$
Furthermore, by the second part of (4.4) and Sobolev's lemma 
$$
\|\div W\|_r\leq \| W_1\|_{r+1}\leq C_r\|\div W\|_r,
\qquad\text{and}\qquad 
\langle W_1\rangle_r\leq C_r \|W_1\|_{r+1}\tag 7.15 
$$
This concludes the proof of Theorem {3.}1.

\subheading{Acknowledgments} I would like to thank Demetrios
Christodoulou, David Ebin and Kate Okikiolu for helpful discussions.

\enddocument